\documentclass[11pt]{amsart}
\usepackage{amscd,a4wide,tikz,amsbsy,mathrsfs,enumerate,mathtools,stmaryrd,tikz-cd}
\usepackage[latin1]{inputenc}
\usepackage{amsmath}
\usepackage{amssymb}
\usepackage{hyperref}
\usepackage{color}
\usepackage{changes}

\definechangesauthor[color=purple]{V}
\definechangesauthor[color=blue]{R}

\usepackage[english]{babel}

\usepackage[normalem]{ulem}

\newtheorem{lemma}{Lemma}[section]
\newtheorem{theorem}[lemma]{Theorem}
\newtheorem{conjecture}[lemma]{Conjecture}
\newtheorem{proposition}[lemma]{Proposition}

\theoremstyle{definition}

\theoremstyle{remark}

\newtheorem{example}[lemma]{Example}
\newtheorem{definition}[lemma]{Definition}

\theoremstyle{plain}

\newcommand{\piphy}{\pi_{\tt{phy}}}
\newcommand{\piint}{\pi_{\tt{int}}}

\newtheorem{myexp}{Example}

\begin{document}

\title[Meyer]{Meyer  sets, Pisot numbers,   and self-similarity  in symbolic   dynamical systems}

\author[V.~Berth\'e]{Val\'erie Berth\'e}
\address{IRIF, Universit\'e Paris Cit\'e - B\^atiment Sophie Germain, Case courrier 7014,
8 Place Aur\'elie Nemours, F-75205 Paris Cedex 13, France}
\email{berthe@irif.fr}
\thanks{This work was supported by the Agence Nationale de la Recherche  through the project  ``Codys''  (ANR-18-CE40-0007), and the
   EPSRC  grant numbers EP/V007459/2 and EP/S010335/1.
 }
\author[R.~Yassawi]{ Reem Yassawi}
\address{School of Mathematical Sciences, Queen Mary University of London, Mile End Road, London, E1 4NS, United Kingdom}
\email{r.yassawi@qmul.ac.uk}


\date{\today}

\begin{abstract}
Aperiodic order  refers to the mathematical formalisation of   quasicrystals.
Substitutions and cut  and project  sets  are   among their main  actors; they also play a key  role in 
  the study of  dynamical systems, whether they   are  symbolic,     generated by tilings, or point sets.
We focus here on  the relations between  quasicrystals   and self-similarity   from an arithmetical and    dynamical  viewpoint,
illustrating   how    efficiently  aperiodic order    irrigates    various domains of mathematics and  theoretical
computer science, on a  journey  from Diophantine  approximation to computability theory. In particular, we see how Pisot numbers allow  the  definition of  simple  model sets, and  how they also intervene  for   scaling factors for   invariance  by multiplication of Meyer sets.
 We  focus in particular  on the characterisation due to Yves  Meyer: any Pisot or Salem number is a  parameter of   dilation that preserves some  Meyer set. 
\end{abstract}

\maketitle


\section{Introduction}

 In his beautiful text   from the book devoted to   the Abel Laureates   \cite{Abel}, Yves Meyer  recalls:  ``I read the  extraordinary book  {\em Ensembles parfaits et S\'{e}ries trigonom\'{e}triques} by J.P. Kahane and R. Salem. I was enthusiastic. The role played by Pisot numbers in the problem of uniqueness for trigonometric series had been discovered by Salem and Zygmund and was detailed in this book. ($\dots$) I was fascinated by the interplay between number theory and harmonic analysis which is so profound and beautiful in this result. I decided to extend the work of Salem and Zygmund to the problem of spectral
synthesis.''  Thus Y. Meyer  began a journey that  would culminate in the  two books \cite{Meyer:70, Meyer:72}, and that would lead  him to  quasicrystals and tilings via model sets.

 We recall that a  real algebraic integer $\beta>1$ is a  \emph{Pisot number}
  (also called  a Pisot-Vijayaraghavan  number)   if all of its other  algebraic conjugates  $\lambda$ satisfy $|\lambda|<1$.
Y. Meyer  thus  started his journey   by  extending  the work of R. Salem and A.  Zygmund to the problem of  {\em spectral synthesis} for    the     Cantor set  $E_\beta$  with $\beta>2$, where
\[E_\beta := \left\{ \sum_i \varepsilon_i \beta ^{-i} \mbox{ with } \varepsilon_ i \in \{0,1\}\right\}.\]  
The aim of  spectral synthesis is to recreate a function from its Fourier coefficients: a   compact set  $\Gamma$  is called a set of spectral synthesis if 
any    continuous bounded function   $f$ with spectrum in $\Gamma$ (i.e.,
whose  Fourier transform has support in $\Gamma$),  can  be approximated in the weak star topology by the set  $\mathcal P_{\Gamma}=\{  \sum_{\gamma \in \Gamma} a_{\gamma} e^{2i  \pi \gamma x} \} $ of   finite trigonometric polynomials with support  in $\Gamma$.
As was similarly done  for  the case of  sets of uniqueness, namely those sets outside of which a Fourier series uniquely determines a function,
Y. Meyer showed that   Pisot numbers also  play a  crucial  role in the problem of spectral synthesis. 
Indeed, for $\beta >2$, Y. Meyer proved   the striking fact  that   $E_\beta$ is a set of spectral synthesis if  $\beta$ is  a Pisot  number. 
Strong spectral synthesis even  holds  with respect to uniform convergence on  compact sets, as stated below.

\begin{theorem}  \cite[Theorem V, Chapter VII]{Meyer:72}  \label{theo:intro} Let $\beta >2$ be a Pisot number.  For each bounded continuous function  $f : {\mathbb R} \rightarrow {\mathbb C}$
whose spectrum lies in $E_{\beta}$,  there exists a sequence 
$(P_k)_k$  of finite trigonometric sums   such that 
\begin{itemize}
\item for each $k$,  the    frequencies  of  $P_k$ (i.e., its spectrum)    belong to the set of finite sums $ \{ \sum_{j=1}^k  \varepsilon_j \beta^{-j} \mbox{ with } \varepsilon_ i \in \{0,1\} \}$,  
\item  the sequence $(P_k)_k$ converges uniformly 
  to $f$ on each compact subset  of ${\mathbb R}$,  and 
  \item     $\sup_{\mathbb R}  |P_k|$ converges to $  \sup_{ \mathbb R} |f|$.
  \end{itemize}    \end{theorem}
  This resolved a longstanding conjecture of R. Salem and A. Zygmund  and extended the result 
that C. Herz   had obtained for the triadic  Cantor set    $E_3$.  
In addition,   the tools that Y.~Meyer   developed allowed him to discern the relationship between the arithmetical nature of a  set   $ \Lambda $, not necessarily compact, and the properties of  the finite trigonometric polynomials with support in $\Lambda$.   As an illustration, let us  start with  the simplest  case $\Lambda = {\mathbb Z}$. The almost periodic  functions with spectrum in $\mathbb Z$, i.e., those that are uniform limits of   trigonometric polynomials in   $\{  \sum_{n \in \mathbb Z} a_{n} e^{2i  \pi n x} \}$, 
 behave like  periodic functions. 
 Y. Meyer introduced a class of sets  $\Lambda$,    considered as   {\em pseudo-lattices},  which have the property that any almost-periodic function with  spectrum (frequencies) in $\Lambda$ behaves like a periodic function.
   Y. Meyer  investigated  these sets $\Lambda$  in the context of Diophantine approximation  for the   study  of  Pisot  and  Salem  numbers in  the books \cite{Meyer:70,Meyer:72}; see also \cite{Meyer:95,Meyer:12}.      
 Here, in view of Theorem \ref{theo:intro},  these sets     are   given by  the discrete grids $\Lambda_{\beta}$  (obtained     from the sets $E_{\beta}$ up to a suitable  rescaling)  with   $$\Lambda_{\beta}=
  \left\{ \sum_{i {\in F}} \varepsilon_i \beta ^{i} \mbox{ with } \varepsilon_ i \in \{0,1\} {, F \mbox{ finite }}\right\};$$ 
 they  provide explicit approximation schemes,   underlying  Theorem \ref{theo:intro}. 
 For specific $\beta$, they  are examples of  sets that  are  now referred to as {\em Meyer} and {\em model} sets.
 Let us briefly   recall      their definitions.

\begin{definition}[Delone and Meyer sets]
A {\em Delone set}  is a   subset   $\Lambda $ of ${\mathbb R}^n$     which  is  both uniformly discrete and relatively dense, i.e.,  there  exist  $r, R > 0$
  so that each ball  with  radius $ r$  contains at most  one point of  $\Lambda$ and each ball of radius $ R$  contains   at least one point of  $\Lambda$.  
A  Delone set $\Lambda $  is    a  {\em Meyer set}
 if  there is a finite set
 $F$ such that 
 $\Lambda - \Lambda \subset \Lambda +F$, or  
 equivalently,  if $\Lambda - \Lambda$  is also a Delone set  \cite{Meyer:72,Lagarias:96,Moody:97}. 
 \end{definition}
 Meyer sets        play   the role  
of   lattices
for  crystalline structures.
With  atoms being  located    at  points  of $\Lambda$, the set $\Lambda - \Lambda$
 corresponds to interatomic distances, and restrictions on interatomic distances  produces    long-range order.

 In  \cite{Meyer:72} Y. Meyer introduced {\em cut and project schemes}, yielding {\em  model sets}, which are explicit constructions of Meyer sets. 
In order to   introduce the idea in the simplest Euclidean setting, we consider  a full rank lattice 
 $L$ in ${\mathbb R}^d$, i.e., a lattice generated by a set of $d$ linearly independent vectors,     together with  a decomposition of 
 ${\mathbb R}^d = {\mathbb R} ^n \times {\mathbb R} ^k $, with $d=n+k$, where  ${\mathbb R} ^n$ is called the {\em physical space} and  ${\mathbb R}^k$ is called the {\em internal space}.
 Let 
   $\piphy$   and $\piint$ denote the natural projections  onto 
${\mathbb R}^n$ and 
${\mathbb R}^k$ respectively.  We further assume that 
 the restriction   $\piphy$  to $L $  is  injective and  that   $\piint(L)$ is dense in ${\mathbb R}^k$.  We project the lattice $L\subset  {\mathbb R}^d$    onto  the  physical space.   The lattice $L$ thus   brings  some order and  regularity  with its  underlying higher-dimensional periodicity; its projection  onto the internal space
provides the  aperiodic structure.

Given a cut and project scheme,
 model sets   are then   formed by projections   together with  a   way   of selecting points.  This selection  is  done  thanks to an
{\em acceptance window} $W$  that lives in the  internal space. 
It  is usually assumed to be  relatively compact  and    with non-empty interior. We now can define Euclidean model sets.
 
 \begin{definition}[Model sets] A subset $\Lambda$ of  ${\mathbb R}^n$
is called  a  \emph{ Euclidean model set}
if there exist a cut and project scheme $({\mathbb R}^n \times {\mathbb R}^k, L)$ 
and  a  relatively compact set  $W$ of $\mathbb R^k$ with non-empty interior such that 
$$\Lambda=\{\piphy(P); \ P \in L ,  \  \piint (P) \in  W \}.$$ 
\end{definition}

Nowadays one allows the internal space ${\mathbb R}^k$ to be replaced by a locally compact Abelian group; this is an important generalisation. For example, the (vertices of the) Penrose tiling is a model set with a non-Euclidean internal space. If we restrict to Euclidean model sets,  the best one can do is to express the Penrose tiling as the union of four translates of Euclidean model sets. Furthermore, the translate of a Euclidean model set is a model set, though not necessarily a Euclidean model set.  

A model set $\Lambda$ is a Delone set: the uniform discreteness comes from the compactness of  the  acceptance window $W$ and the relative denseness from its non-empty interior. One  of the main  features of model sets is  that, under simple extra restrictions, they      provide    pure  point  diffraction (or, in more dynamical  terms, pure discrete spectrum),  such as  highlighted in \cite{Moody:20}.

Y. Meyer established the following powerful  connections. For details, see   \cite[Theorem IV  Chapter 2, p.48]{Meyer:72}, \cite{Meyer:95}
    and for more about 
model  sets,  see for instance  the surveys 
\cite{Moody:97,Moody:00,Lee&Moody:2001} and \cite[Chapter 7]{AperBook}. 
 \begin{theorem}
   Model sets are Meyer sets. Conversely,  if $\Lambda$  is a Meyer set,  then it is a  relatively dense   subset  of finitely many translates of a Euclidean model set, i.e., there exist a Euclidean model set $M$  
  and a finite set $F$ such that $ \Lambda \subset M+F$. 
  \end{theorem} 

 The work of Y. Meyer on model sets found  a second life  with the discovery of  quasicrystals.  In 1974, R. Penrose discovered his tilings\footnote{A {\em tiling}  is a  covering of  a surface  (usually  a Euclidean plane)  using  copies of geometric tiles, placed next to each other, without holes or overlaps.},
and in 1981 the connection between these tilings and Y. Meyer's work was made explicit  thanks to N. de Bruijn  who proved  in \cite{debruijn} that the vertices in the Penrose tilings formed a model set. 
  M.  Duneau, D. Gratias, A.  Katz, and R. V.  Moody then discovered that the {\em quasicrystals} described by D. Shechtman can be modeled using model sets.

A  quasicrystal is a physical solid whose atoms are arranged as a  translation invariant lattice which does  not satisfy the crystallographic restriction.  D. Shechtman's discovery of quasicrystals 
 in 1984  \cite{QC}    mobilised    the 
    material sciences and physics  to  provide  suitable models for the description of quasicrystals. 
    Mathematicians came to realise that the theory of aperiodic tilings, with roots in both the theory of Wang tiles (Section \ref{sec:comput}), and also
     in  recreational mathematics,   
where  R. Penrose's aperiodic    tilings displayed  pentagonal symmetry and   thus could not satisfy the crystallographic restriction, provided a context in which to study mathematical models of quasicrystals. In particular, one studies  these physical objects by looking at an appropriate cut and project set of their atom locations, and this leaves us with Meyer sets.

Nowadays the term {\em aperiodic order}, which first appeared in  Moody's work \cite{Moody:97}, refers to the study of mathematical models of quasicrystals, continuous or discrete.  As with the study of quasicrystals, aperiodic order  is the study of mathematical objects which simultaneously have
   long-range order
and  aperiodicity.   The connections to Y. Meyer's   concepts and  to aperiodic  tilings   are now well recognised and  established, in particular thanks to   R. V. Moody \cite{Moody:00} and to J. Lagarias \cite{Lagarias:96}, and the  fundamental contributions to aperiodic order   brought by Y. Meyer  developed  in his founding books  \cite{Meyer:70,Meyer:72}   continue to nourish the field. See also \cite{GQ:20,Moody:20,OU:20} in this volume. This area has   now   exploded and covers a  wide range of domains  such as illustrated by the  books \cite{sene,CRM:2000,MAP:15} and the   book series {\em Aperiodic order} by  M. Baake and U. Grimm, in particular  \cite{AperBook}.

The notion of  long-range order is crucial here:  crystallographers  now  mean  an  ordering of atoms   which  produces a diffraction pattern with sharp bright spots, such as 
confirmed by  the general definition of crystals adopted by  the  Crystallographic Union  in 1992.
These diffraction patterns translate mathematically to the study of the spectrum of point sets.
 In a nutshell,  aperiodic order  goes with 
diffraction and pure discrete spectrum, as   highlighted  in   the contribution by R. V. Moody  in this volume \cite{Moody:20}. But the    beauty   and the strength of  aperiodic order is that it   goes well beyond  spectral considerations.

There is  in particular  a  further crucial  aspect here,   namely   the prevalence of    self-similarity in aperiodic order. In mathematics, a self-similar object is exactly  or   approximately similar to a part of itself, or it can be divided into smaller copies of itself; 
 hence  it looks roughly the same on any scale and  is invariant upon being scaled larger or smaller. It can be that parts of a figure are small replicas of the whole, 
 or else a self-similar group.   The paradigmatic example is a  compact set $K$ which is the  solution of  an  iterated function system, i.e., 
  $K= \cup S_i(K)$,  where  $i$ belongs to  a finite set  and  where the maps $S_i$ are contractions.
A classical way to generate self-similar objects is by iterating a rule.  For example, one can iterate functions, to obtain fractals such as the Sierpinski gasket, one can substitute a side of a triangle by a polygonal line, such as in the Koch curve.  Among the most classical   methods for  generating mathematical  models of quasicrystals
are  {\em substitutions}  which will be  the main  object of Section \ref{sec:subs}.
In aperiodic  order we iterate these substitutions to obtain self-similar sequences, points sets or tilings.   Indeed, in  tiling theory, a tile represents an atom  and a tiling represents  large  atomic configurations.
 In point sets,  points  are  idealisations of atomic positions,    and    letters  are used to code     interatomic  distances, and, considering specially chosen control points in the tiles  allows one to  equivalently  associate Delone   point sets with tilings.

Algebraicity  provides a particularly relevant  arithmetic expression of self-similarity, and algebraic numbers, and  particularly   Pisot numbers, are  at the   heart of aperiodic order.   The 	aim of this survey is twofold.  We   want to  stress the impact of  Pisot  numbers in   aperiodic  order,   and also   to highlight the explosion of   the involved domains,  by focusing on an approach based on dynamical systems.   The involved mathematical fields    cover here  an impressive spectrum, with   the  toolbox of   aperiodic order  containing  
 Diophantine  approximation and  arithmetic, together with  topological methods via topological  dynamics or cohomology, computabillity theory, 
spectral theory, algebra, discrete geometry and harmonic analysis.
There are several  reasons that explain this abundance.      Firstly,  aperiodic order  applies    to  various types of    quasiperiodic structures which  can  serve as models for quasicrystalline structures,    whether they are  geometric with 
 tilings   and point sets,  or    symbolic   with   infinite  words.  
Secondly,    the  two main  types of  mathematical models  for    quasicrystals,  namely 
 model sets  obtained via   cut and project schemes, and  substitutions,    are  prominent  concepts 
 in the   study of dynamical systems,   allowing dynamicists to express  self-similarity in  both  arithmetic  and    geometric  terms.  As an illustration, consider the  role played by
 substitutions    
 via   renormalization  in the study of interval exchanges \cite{yoccoz}.

\bigskip

 In  this survey we made  the choice  
  to focus  on  several   complementary    topics that illustrate the   strength of the concepts of   cut and project schemes and substitutions in  the dynamics of aperiodic order. 
 In Section \ref{sec:Pisot} we see how Pisot numbers allow the definition of simple model sets, and   we recall   the following characterisation, due to Y. Meyer. Any Pisot or Salem number is the parameter of  dilation preserving some Meyer set: this  
is  Theorem \ref{thm:meyer} and it will serve as a common thread.   In Section  \ref{numeration} we discuss how  Meyer  and model  sets  are provided by  Pisot numbers  and  beta-numeration.
 Substitutions, their definitions,  stability by multiplication and self-similarity  are   considered in Section \ref{sec:subs}.   In Section \ref{sec:Pisot-Sturmian-Rauzy-Sadic} we discuss  the  spectral properties  of substitutions and the Pisot conjecture, simple one-dimensional model sets given by 
  Sturmian words and more exotic acceptance windows provided by Rauzy fractals.
 In Section \ref{sec:comput}, we discuss how  computabillity enters the picture    with the question of  the existence of   local rules for aperiodic tilings  and   the domino problem.
 
 \bigskip

 \subsection*{Acknowledgements}
 We would like to thank  Nathalie Aubrun, Thomas Fernique, Pierre-Antoine Guih\'eneuf, \'Etienne Moutot,  Wolfgang Steiner, Nicolae Strungaru  and Pascal Vanier   for  their careful reading and their  suggestions. We also thank Uwe Grimm for kindly  allowing us to use Figure \ref{Uwe}.  Lastly, we thank   the  anonymous referees for their careful reading and their constructive suggestions.

\section{Pisot numbers and aperiodic order  }\label{sec:Pisot}

Algebraic numbers, and  particularly   Pisot numbers, are  at the   heart of aperiodic order. Algebraicity intervenes in two ways: firstly, Pisot numbers  allow us to define  model  sets  $\Lambda$, and secondly,  {\em scaling invariance}, i.e.,  the existence of a scaling factor  $\beta$ with $\beta \Lambda\subset \Lambda$, is a  step towards self-similarity.
This is seen in  Theorem \ref{thm:meyer}, a striking result by Y. Meyer  which has  led  to a rich theory   in the study of associated dynamical systems, some of which  we will exposit in later sections.
Let us  first recall  the following definition. 

\begin{definition}[Salem number]
A real algebraic integer $\beta>1$ is a 
\emph{Salem number}  if  at least  one of its algebraic  conjugates $\lambda$ satisfies $|\lambda|=1$ while its  other  conjugates  have modulus  smaller than $1$.
\end{definition}

The following result will be   a  common thread  throughout  this survey. See also \cite{Lagarias:99} for a similar statement for Delone sets.

\begin{theorem} \cite[Theorem 6]{Meyer:95}  \label{thm:meyer}
 If   $\Lambda$ is a Meyer set,  $\beta >1$ is a real number and if $\beta \Lambda \subset \Lambda$, then $\beta$ is either   a   Pisot  number
or a Salem number.   Conversely, for  each dimension $n$  and for  each Pisot or   Salem number $\beta$, there  exists  a model set  $\Lambda \subset {\mathbb R}^n$ such
that $\beta \Lambda \subset \Lambda$. 
\end{theorem}

The proof of the first statement    goes via the concept of    {\em harmonious sets} which  provide  a further   characterisation by duality  of  Meyer sets.

\begin{definition}[Harmonious set]
A harmonious set   is a subset $F$ of a locally compact abelian group $G$ such that  every weak character\footnote{A weak or algebraic character  is a homomorphism  for which no continuity is required.}    on the subgroup generated by $F$
  may be approximated uniformly by a continuous character on $G$. 
  \end{definition}

For more on the subject, see \cite{Meyer:72}, \cite{Meyer:95},  \cite{Meyer:20}. 
 One considers here  characters of the additive group generated by the Meyer set $\Lambda$, equipped with the discrete topology.   The following beautiful result tells us that harmonious sets are non-other than Meyer sets (see     \cite[Chap. II Section 5]{Meyer:72},  \cite[Theorem 4]{Meyer:95}).
  \begin{theorem} \label{thm:harmo1}
Let $ \Lambda$ be a Delone set.   Then $\Lambda $ is  a Meyer set if and only if $ \Lambda$ is harmonious.
\end{theorem}

Using Theorem \ref{thm:harmo1},
we can prove the first part of Theorem \ref{thm:meyer}.  We follow here \cite{Meyer:95}. Take a Meyer set $\Lambda$  and let $\beta >1$   be such that    $\beta \Lambda \subset \Lambda$. Take a non-zero  element $\lambda_0$ in $\Lambda$. Then the set  $\{\beta^k  \lambda_0: k \geq 0\} $ is included in $ \Lambda $.
We then  note   that  a subset of a  harmonious set  is again a harmonious set,  and  that      the set  $\{\beta^k  \lambda_0: k \geq 0\}$ is harmonious in  the one-dimensional space $\lambda_0 {\mathbb R}$. 
Since  $\Lambda$ is  a Meyer  set, then   $\Lambda$ is  a harmonious set.
We now can conclude by using the following.
 \begin{theorem} \cite[Theorem XX, Chapter 1]{Meyer:72} \label{theo:harmo}
 Let $\beta>1$.   The set  of powers  $\{\beta^k,  \ k \geq 0\}$ is harmonious if and only if $\beta$ is a  Pisot  or a Salem number. 
 \end{theorem} 

This sketch of   proof is an illustration of the  fact that ``weak characters pave the road which goes from S. Bochner to D. Shechtman'', as 
nicely phrased  by Y. Meyer in  \cite{Meyer:20}. He moreover adds: ``On the way we are visiting  {\em coherent sets of frequencies,  harmonious sets},(...)  and finally we arrive at {\em model sets}'', a visit  to coherent sets  of  frequencies that we will also  make   in Section \ref{numeration}.
We will not say more on the first statement of Theorem \ref{thm:meyer}, but  rather concentrate on the second one, that is, the construction of the model set related to $\beta$.
  But before providing some   elements of the  proof of the second statement,  we  first  illustrate it  with a paradigmatic and  illustrative  example of how 
 algebraic numbers  provide us with  the first and simplest examples of model  sets.  We later revisit this running example in other settings where aperiodic order is studied. Although it has a  simplicity   which does not   reflect the subtleties   needed in  the general case,  it nevertheless   gives  the reader a flavour of why these results  may be true. See also \cite{GQ:20} in this volume.

 \begin{example}\label{ex:Fibo}
 Consider  the golden ratio $\varphi$, which is the largest  root of $x^2-x-1=0$, and its algebraic conjugate $\varphi'$, and define the lattice \[\mathcal L:=\{ m(1,1)+n(\varphi, \varphi'): m,n \in \mathbb Z\}=\{ (m+n\varphi, m+n\varphi'): m,n \in \mathbb Z\} .    \] In the first coordinate (the physical space), the lattice $\mathcal L$  projects injectively  by  $\piphy(\mathcal L)$ onto $\mathbb Z[\varphi]$ in $\mathbb R$, and in the second (the internal space), its projection $\piint(\mathcal L)$ is dense in ${\mathbb R}$.
Now take  $W=[-1, \varphi-1]$; then 
\[ \Lambda:= \{   m+n\varphi:  m,n \in \mathbb Z, -1 \leq m+n\varphi'\leq \varphi -1  \}\]
is a model set; see Figure \ref{Uwe}.

 Moreover, this  model   set    is scaling invariant under   multiplication   by $\varphi$. We use the fact that  $\varphi$   is a  quadratic unit.
 Indeed,   assume that  $m+n \varphi  \in \Lambda$. 
As $\varphi^2= \varphi+1$, we have   $ \varphi (m +n \varphi)=n+  (m+n)  \varphi$. Similarly, 
  $n+  (m+n)  \varphi'= \varphi' (m +n \varphi')$. Multiplying the given inequality
   $-1 \leq m+n\varphi'\leq \varphi -1 $ by $\varphi'$  and using $\varphi \varphi'=-1$ and  $\varphi+\varphi'=1$
   yields  $\varphi -1 = - \varphi' \geq n+  (m+n)  \varphi' \geq -1 - \varphi' \geq -1$, 
 proving the   stability of $\Lambda$  under  multiplication by $\varphi$.

Finally, it can be shown  that this cut and project scheme generates both a one-dimensional tiling of $\mathbb R$ using two distinct tiles, and if we label these tiles $A$ and $B$, we obtain  an {\em infinite word}, i.e., an infinite concatenation of the letters $A$ and $B$.  We will see in Section \ref{sec:subs} that this discrete object is the fixed point of a {\em substitution}, and a discrete model of self-similar aperiodicity.  See also \cite{GQ:20,OU:20} in this volume.

\begin{figure}\label{Uwe}
	\includegraphics[width=.5\textwidth]{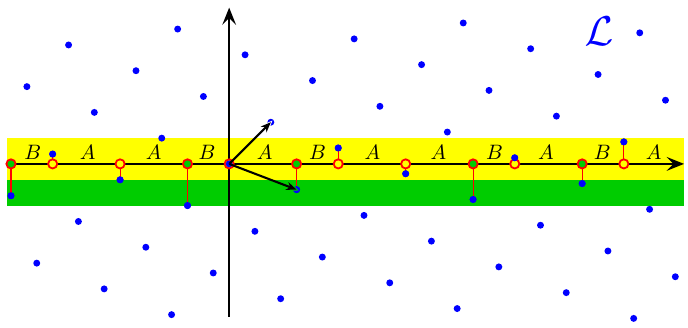}
	\caption{ A model set generated by the golden mean $\varphi$; it  consists of the circled points on the horizontal axis.  The tile length of $A$ is $\varphi$ and that of $B$ is $1$.}
\end{figure}

\end{example}

More generally,  given two distinct irrational parameters $\nu, \epsilon$  (which   are not necessarily  algebraic)  and  a bounded acceptance interval  $W$, one similarly  gets a  one-dimensional  model set
   $\{a+b  \nu  \mid a,b \in {\mathbb Z}, a+b \epsilon \in W\} \subset{\mathbb Z}[\nu]$. As in the example above, such    cut and project schemes,  obtained by selecting points with integer coordinates located within a bounded distance of a line, and  thus  having
   one-dimensional physical and 
    internal  spaces,  lead    to one-dimensional  tilings.
Depending on whether or not the length of the window belongs
   to ${\mathbb Z} [\epsilon]$,
   the lengths  in the tiling take two or three values, according to the size of the   acceptance window  \cite{GMP:03},
    and  if  we code the intervals with distinct letters, we obtain  two- or   three-letter   infinite words.
  Such   codings  correspond either to Sturmian words (see Section \ref{subsec:sturm}), in the two interval case,  or to codings of so-called three-interval exchanges in the three-interval case.
Special attention has been  given to the combinatorial properties   of this   particularly  rich class of   infinite  words: see for instance   \cite{GMP:03,GMP:06} and the survey \cite{PelMas:07}.
Such a  model set    has    a priori no  scale invariance. 
Imposing the requirement of self-similarity implies
that  $\epsilon$   is a quadratic integer and $ \nu= - \epsilon'$, i.e., the algebraic conjugate of $\epsilon$  and the scaling factor $\gamma$
must be a quadratic Pisot number in the same algebraic field $\mathbb{Q}(\epsilon)$ (see \cite{GMP:03}).  This  brings us back to Y. Meyer's connection between self-similar  model sets and  Pisot numbers  in Theorem \ref{thm:meyer}.

Pisot  numbers  allow the extension   of  such a   simple   and seminal  construction to  higher dimensions,  which   is   a key point in
the proof of the second statement of Theorem \ref{thm:meyer},  which  we now summarise 
 from \cite{Meyer:95}, and  which is  illustrated with Example  \ref{ex:Fibo}.
 Briefly,
when  $\beta$ is  a Pisot  or  a Salem number, the     construction of   a model set is algebraically realised  by taking the lattice $\mathcal L$
 to be the   ring of algebraic integers in 
${\mathbb Q}( \beta)$ and   the 
 canonical embeddings  $\sigma_i$ associated with   the (other)  algebraic conjugates   of $\beta$  provide the projection $\piint$ 
onto the internal  space, with  the projection $\piphy$ onto the physical space    being   the identity.   The  acceptance window is  given
by   the conditions $|\sigma_i(x)| \leq 1$.
The model set is thus the  set of algebraic integers  $\lambda$ in ${\mathbb Q}( \beta)$ such that   $|\sigma_i(\lambda)| \leq 1$
for all  $i \geq 2$.  
This   approach  by Y. Meyer  is now  classical and  has  had a   particularly  fruitful influence on    subsequent constructions.

\section{Pisot numbers and  beta-numeration}\label{numeration}
 Arithmetic  examples of  sets which clearly have scaling invariance  are provided by   the sets  $\Lambda_{\beta}=  \left\{ \sum_{i {\in F}} \varepsilon_i \beta ^{i} \mbox{ with } \varepsilon_ i \in \{0,1\} {, F \mbox{ finite}}\right\}$
  that we  met in the introduction.   These grids have  remarkable  properties that  lead to
 the notion of  a {\em coherent set of frequencies}  (see \cite[p.110]{Meyer:72}).

\begin{definition}[Coherent set of frequencies]
A  set  $\Lambda$ of real numbers   is    a {\em coherent set of frequencies}   if  there  exist a   constant $C>0$ and a compact set $K$ 
of real numbers such that 
$$ \sup_{\mathbb R}  | P(x) | \leq  C \sup_K |P(x) | $$
  for all  finite trigonometric sums $P$ with frequencies  in  $\Lambda$.
\end{definition}
 For more on  the contributions of  Y. Meyer in this setting, see  the enlightening  text by A.~Cohen in \cite{Abel} and see  also  \cite{Meyer:20} which   revisits    the relationship between   coherent sets  of frequencies and Bochner's  property.

The following result should be compared with Theorem \ref{theo:harmo}.
\begin{theorem} \cite[Theorem IV, Chapter IV; Proposition 7, Section 8, Chapter II]{Meyer:72} \label{thm:harmo}
{Let $\beta>1$.}
The following are equivalent:
\begin{itemize}
\item   the number $\beta$ is a Pisot number;
\item 
the set $ \Lambda_{\beta}$ is harmonious;
\item
the set $ \Lambda_{\beta}$ is  a coherent set of frequencies.
\end{itemize}
\end{theorem}

Note that this statement  does not mean that $\Lambda_\beta$ is a Meyer set, since it is not relatively  dense  when $\beta >2$, and hence not Delone.
 Indeed, consider   the strictly increasing  sequence $(x_n)$ of  numbers that have at least one representation as an element of $\Lambda_\beta$, and 
then 
its {\em first difference sequence} $(x_{n+1}-x_n)$.   The following holds.

\begin{proposition}\cite[Lemma 1.1]{FW:02}\label{prop:firstdiff}
If $\beta >2$, and $(x_n)$ is an increasing sequence of  numbers that have at least one representation as an element of $\Lambda_\beta$, then its first difference sequence    takes arbitrarily  large values.
\end{proposition}

 More generally, sets 
  of  finite  sums of the form   $\Lambda_F=\{\sum_i\varepsilon_i \beta ^i \mid \varepsilon_ i \in F\}$,   for specific  finite sets  $F$ of  non-negative integers,  when $\beta$ is assumed to be a Pisot number, have been widely considered.
Researchers in  arithmetic   and numeration  came indeed to  notions similar to Meyer sets independently via base expansions.
   In particular,  the    study of  the  first difference  sequence was initiated in the fundamental work    \cite{EJK:90} for $ \Lambda_F=\Lambda_{\beta}$  with    $1 < \beta <2$ and $F=\{0,1\}$.   The study of the  limit inferior  of the   first difference sequence has  attracted  much  attention.
  In particular, $\liminf x_{n+1}-x_n>0$ 
  if   and only if  $\beta $ is  a Pisot number; see   \cite{Bugeaud:96} and also \cite{EJK:98}.

       Returning to  the golden mean  Example \ref{ex:Fibo},  if we take $\beta=\varphi$ and $F=\{0,1 \}$, one can check that $\Lambda_F = \mathbb Z[\varphi]$ and that the entries of the   first difference sequence are either 1 or $\varphi -1$. Thus $\Lambda_F$ is nothing other than our lattice $\mathcal L$ and the entries $1$ and $\varphi-1$ of the difference sequence give us the lengths of the tiles that we saw in Figure \ref{Uwe},  when scaled by a factor of $\varphi$.

   The  case where $\beta \geq 2$    and    $F= \{1,\cdots, m\}$  has  also  been widely studied.
See   e.g.  \cite{AK:13}, and the references therein,   for the study of 
     set of all differences       $\sum  \varepsilon_i \beta ^i$ where $\varepsilon_i \in \{0, 1 , \pm 1, \cdots, \pm m\}$   in relation with singular Bernoulli convolution distributions;    note that  this latter set  gives information about $\Lambda_F- \Lambda_F$; 
in particular the   absence of accumulation points   in this set    is  proved to be equivalent  to  $ \beta$  being  Pisot  or  $\beta \geq m+1$  \cite[Theorem 1.1]{AK:13}.
In  the  language of  aperiodic order,  this gives  the following:  
 \begin{proposition} \label{prop:F} Let $\beta>1$.
 For  $F=\{0, 1 , \cdots, m\}$,   the  set  $ \{P(\beta); \ P \in  F[X]\}$   is   a Meyer set if and only if  $\beta  \leq m+1$ and $\beta$ is a Pisot number.
 \end{proposition}

 Let us provide a  sketch of the proof of Proposition \ref{prop:F}. We first  briefly   recall how   the Pisot property intervenes here as   being  crucial in guaranteeing the uniform discreteness of $ \Lambda_F= \{P(\beta); \ P \in  F[X]\}$ for $F$   a finite set of  integers.  This is a  classical argument  known as  Garsia's Lemma; see   \cite[Lemma 1.51]{Garsia} or \cite[Lemma 6.6]{Sol97}.
If $\beta$ is assumed to be Pisot,
there exists $C$ such that $|P_1(\beta ^{(i)} ) -P_2(\beta ^{(i)})| \leq C$
for any algebraic conjugate $\beta^{(i)} \neq \beta$ and any  $P_1, P_2 \in F[X]$. Next,
since the integer $\prod_{i}(P_1-P_2)(\beta ^{(i)})  (P_1-P_2)(\beta) $
 is non-zero  for $P_1, P_2 \in F[X]$ with $P_1(\beta) \neq P_2(\beta) $, we deduce
a positive uniform lower bound    for $P_1(\beta ) -P_2(\beta)$
with $P_1(\beta) \neq P_2 (\beta)$, from which we conclude uniform discreteness.
 
It remains to consider the property of  relative denseness where the cases  $\beta\leq m+1$ and $\beta >m+1$ have to be distinguished  by \cite[Lemma 1.1]{FW:02}.
Indeed, if   $\beta   > m+1$,  the difference sequence  $(x_{n+1}-x_n)$ takes arbitrarily  large values  (this extends Proposition  \ref{prop:firstdiff}), which prevents relative denseness (see the proof of \cite[Lemma 2.1]{EK} such as stressed in   \cite{FW:02}). 
    Next if  $\beta \leq m+1$ and 
$\beta$  is  a  Pisot number,
    the difference sequence  can take only finitely many distinct values;   it   even  has a strong  combinatorial structure such as described in    \cite{FW:02}, where it is shown that the
difference sequence can   be generated by a substitution over a finite alphabet (see Section \ref{subsec:subsdef}).   

Lastly,   assume  that  $ \Lambda_F= \{P(\beta); \ P \in  F[X]\}$ is a Meyer set.   As  $\beta   \Lambda_F \subset  \Lambda_F$, we deduce from Theorem  \ref{thm:meyer}
that $\beta $ is a Pisot or a Salem number.   By the same proof as that of Theorem \ref{thm:meyer} and also  by Theorem \ref{thm:harmo1},
  $ \Lambda_F$ is harmonious,  and so $\Lambda_{\beta}$ is  also harmonious, from which we deduce that $\beta$ is a  Pisot number, by  Theorem \ref{thm:harmo}.

\bigskip

Thus researchers  have   studied the sets $\Lambda_F$   for over fifty years,  due to their connections to quasi-crystals,  infinite Bernoulli convolutions 
and expansions in non-integer bases.  Dynamicists  have also studied these objects, via the  {\em beta-transformation}.
The  {\em beta-numeration}  is a numeration system that allows the  greedy representation   of    real numbers 
as  sums of powers  of  a real number $\beta$, with $\beta >1$,  in the same way as  real numbers  can be   expanded in
decimal  numeration.  
It was introduced  and  studied  by  A. Renyi  \cite{Renyi} and W. Parry  \cite{Parry} in dynamical terms.
Indeed, given  a   positive  irrational  number $\beta>1$,   consider the  $\beta$-transformation $T_{\beta} 
$ acting on the unit interval defined by $T_{\beta}(x)= \{ \beta x\}$.     Digits   of $x$ in its $\beta$-expansion are  then   obtained as $a_i= \lfloor \beta  T_\beta^{-i}x \rfloor$, for $i  \geq 1$ (hence $a_i \leq \lfloor \beta\rfloor$), which yields
$x= \sum_{i \geq 1}  a_i  \beta^{-i}$. 
The beta-numeration has  been particularly well  studied  when $\beta$ is assumed  to be a Pisot number,  both from a dynamical
and a combinatorial viewpoint; there are  indeed specific properties that occur when $\beta$ is a Pisot number.  Returning  to Example \ref{ex:Fibo}  where  $\beta$ is the golden mean,  digits are in $\{0,1\}$ and   the product of two consecutive digits  in the   $\beta$-expansion  equals 0 (since we use the greedy representation). For Pisot numbers $\beta$, the appropriate generalisation of this statement is that the   finite sequences of  digits   that occur  in  the possible representations have  a  simple expression in terms of finite automata (see e.g. \cite{Berstel&Seebold:2002}). From a dynamical viewpoint, 
points in the interval are replaced by the sequence of coefficients in their beta-expansion.

Moving  from the unit interval to $\mathbb R ^+$, one can also  expand non-negative  real numbers.
In this context,  the    countable   set of {\em beta-integers}
${\mathbb Z}_{\beta}$ 
  plays a particular role: elements of ${\mathbb Z}_{\beta}$  are the 
integer parts  of beta-expansions, that is, they  are  polynomials in  
$\beta$  with coefficients  in $ \{0, \cdots , \lfloor \beta \rfloor \}$  with  combinatorial  constraints on their digits determined by the algebraicity of $\beta$.
 In this  dynamical setting,   the constraints  are given on    the digits $\varepsilon_i$  that can appear when expressing the beta-expansion of a  number. These constraints are imposed to guarantee uniqueness of expansions, and  are particularly easy to express when $\beta$ is an algebraic integer.
For example, returning to  Example \ref{ex:Fibo}, the equation $\varphi^2=\varphi+1$ tells us that a fragment $\varphi^{n+1}+\varphi^n$
 of  $\varphi$-expansion  can be replaced by  $\varphi^{n+2}$. For this example,  we consider ${\mathbb Z}_{\varphi}$ to be polynomials in  $\varphi$  with coefficients  in $ \{0, 1 \}$, but where we do not see consecutive monomials. 
In the case where $1<\beta <2$, one sometimes has   ${\mathbb Z}_{\beta} = \Lambda_{ \beta }$,  for example in the case of  our running example  where  $\beta$   is equal to the golden ratio. However in general one only has ${\mathbb Z}_{\beta} \subset \Lambda_{ \beta }$; see \cite{FrougnySolomyak92}.

This set ${\mathbb Z}_{\beta}$  is  self-similar,  in that $\beta {\mathbb Z}_{\beta} \subset  {\mathbb Z}_{\beta}$.
As    stressed  by Y. Meyer in \cite{Meyer:12}, ``Model sets are non-uniform grids'',  and beta-numeration    has   revealed  itself as a
   very efficient tool
for the modeling of families of  quasicrystals thanks to   the beta-grids  $\sum  {\mathbb Z}_{\beta}  {\mathbf e}_i$, where 
the  ${\mathbf e}_i$ are  vectors  of the canonical basis (see  e.g.
\cite{BFGK:98,verger} and  also   \cite{Thurston}).
The next statement, due to  \cite{BFGK:98,verger}, and which also follows from arguments in \cite{Schmidt:80} and  \cite[Proposition 2]{FrougnySolomyak92},   is  again reminiscent of Theorem \ref{thm:meyer}.

\begin{proposition} 
 If $\beta$ is a Pisot number,  then ${\mathbb Z}_{\beta}$ is   uniformly discrete, and 
$  {\mathbb Z} _{\beta } \cup (- {\mathbb Z} _{\beta })$
is a Meyer set. Conversely,  if  $  {\mathbb Z} _{\beta } \cup (- {\mathbb Z} _{\beta })$ is a Meyer set,
then 
$\beta $   has to be  a Pisot or a Salem number.
 \end{proposition}

 Furthermore, the    set    ${\mathbb Z}_{\beta}$ can be endowed  with 
laws that    are close to multiplication and addition. This allows one to recover properties that  have the flavour of 
the  algebraic rules of lattices. For  some families of  Pisot  numbers  $\beta$, mainly Pisot quadratic units and for  some cubic Pisot numbers,  an internal law  can even be produced
formalising this quasi-stability under subtraction and  multiplication:     the set of beta-integers  is then  endowed with an addition law and a multiplication law, which are compatible with the combinatorial structure
which   makes it a  quasiring    \cite{BFGK:98}.
   Let us quote also  \cite{ABF:05},  where   a method based on automata is given to find a set $F$ such that
 ${\mathbb Z}_{\beta} - {\mathbb Z}_{\beta} \subset  {\mathbb Z}_{\beta}+F$.

\section{Substitutive dynamics  and self-similarity} \label{sec:subs}

We have seen  in Section \ref{sec:Pisot} that infinite words  occur naturally as codings of  one-dimensional quasicrystals.  There  is a particularly   simple  yet 
 rich combinatorial  object  that   allows  the expression   of 
scale  invariance and  self-similarity properties:  this is  the notion of  a  substitution.
We define a substitution rigorously below, but  roughly speaking,  a  {\em substitution}  is  an inflation  rule, either  combinatorial or  geometric,    that, following the cases,    replaces  a letter by a  word, or  a  tile by a geometric pattern made of a  finite union of tiles.
Iteration of  substitutions generates   hierarchically ordered  structures   such as  infinite words,  shifts, point sets, or else  tilings, that all display strong self-similarity properties.

 Substitutions originated  with the first works of A. Thue   in 1906, with the so-called  Prouhet-Thue-Morse substitution, which is an illustrious example of an  automatic sequence \cite{ASBook}. Substitutions generate symbolic dynamical systems and symbolic dynamics is  the setting in which we define them.

   Symbolic dynamical systems are discrete dynamical systems where the space is a Cantor space. Typical symbolic systems are {\em shifts} where the space $X$ consists of infinite words with entries  from a finite alphabet, and the dynamics is  the shift map  $S$  moving entries of an infinite  word one unit to the left.   The set $X$ must be shift-invariant  and  closed in  the subspace topology  inherited from the product topology on $\mathcal A^{\mathbb Z}$, where $\mathcal A$ is given the discrete topology.
    Shifts provide representations of  dynamical systems: given a dynamical system we can discretise the state space using a finite   partition, and code trajectories  of points as sequences of symbols, as first done   by J. Hadamard   \cite{Hadamard} and M. Morse  \cite{Morse}. The Jewett-Krieger theorem gives us very general conditions in which we can represent a dynamical system as a shift.

  Substitutions  soon  outreached in an  unexpected way  this  combinatorial    and dynamical   setting.      Indeed, beginning in the  sixties,  substitutions   were 
 used
 to    input computation in tilings, enabling researchers to produce the first examples of    aperiodic  tilings,   such as   Robinson
 tilings,  and thereby  proving the undecidability of  the domino problem (see  Section \ref{sec:comput}).
Also, in the eighties, {\em Pisot  substitutions}     attracted  much attention in the context of  aperiodic order  
   as  mathematical models  of   quasicrystals  displaying  self-similarity,    yielding  in particular   the Pisot  substitution  conjecture (see Section \ref{subsec:Pisot}).

\subsection{Substitutions: how dynamicists came to model sets}
\label{subsec:subsdef} 

Let $\mathcal A$ be a finite alphabet and let ${\mathcal  A}^+$ be the set  of 
non-empty finite words on ${\mathcal  A}$. 
A
(word) {\em substitution} is a map $\sigma: {\mathcal  A} \rightarrow{\mathcal  A}^+$ 
which extends to   $\sigma: {\mathcal  A}^+ \rightarrow{\mathcal  A}^+$ and  $\sigma: {\mathcal  A}^\mathbb N \rightarrow{\mathcal  A}^\mathbb N$
 via concatenation.  A substitution defines a {\em language}, which is the set of all finite words $w_0\dots w_k$ which appear as a subword of   $\sigma^n(a)$, for $a\in \mathcal A$ and $n\in \mathbb N$. The  shift space $X_{\sigma}$    is   the set of infinite words  all of whose finite subwords belong to the  language of $\sigma$, and with the shift map $S:X_\sigma\rightarrow X_\sigma$, the pair $(X_{\sigma},S)$ is called a (substitution) shift.

  \begin{example} \label{ex:Fibo2}
 As a   prominent   example  of  a word  substitution in aperiodic order,  consider the   {\em  Fibonacci substitution} $\sigma$  defined  on the alphabet  ${\mathcal A} = \{A,B\}$ as  $\sigma(A)=AB$ and $\sigma(B)=A$. 
The sequence of finite words $(\sigma^n(A))_n$ are nested, i.e.,
 $\sigma^n(A)$ is a prefix of $\sigma^{n+1}(A)$  for each $n$, and  it converges to    the so-called {\em Fibonacci  word},  whose first terms are  
\[
ABAABABAABAABABAABABAABAABABAABAABABAABABAABAABABA\ldots
\]

  The Fibonacci word    is   a one-dimensional  discrete model of a  quasicrystal; indeed,  one of the first  models that have been considered.  For, it   codes   the one-dimensional tiling obtained   as  a model set  via  the
  selection  strip in ${\mathbb Z}^2$   illustrated in Figure \ref{Uwe}, noting that the Fibonacci word will be the coding of the right part of the tiling seen there, starting at the origin;  see below for how to generate the actual tiling.  The Fibonacci word    also  corresponds to a  discretization   of a line in discrete geometry \cite{Rosenfeld}, and  from  a symbolic dynamics viewpoint, 
 it belongs to the  family of Sturmian words which are discussed in Section \ref{subsec:sturm}.  For more on the  Fibonacci word, see e.g. \cite{Lothaire:2002,PytheasFogg:2002}.
 \end{example}

The {\em incidence matrix} $M_{\sigma} $ of  a substitution $\sigma$   is   the square matrix  whose entries  count the number of occurrences of letters in
 the  images of letters. For example, for  the Fibonacci substitution $\sigma$, we have
$M_{\sigma}= \left[\begin{array}{ll}
1& 1 \\
1 & 0
\end{array}\right]$.
This matrix, called the {\em abelianization} of   $\sigma$ as the order of letters in substitution words is forgotten, allows us to  
apply tools from linear algebra to analyse the self-similarity properties of the  Fibonacci word, as well as its combinatorial  and dynamical  properties.
The Perron--Frobenius theorem tells us that if  $M_\sigma$ is primitive, then it
 admits a strictly dominant eigenvalue   that is  positive,  and 
 this  eigenvalue  is called   the  {\em inflation factor}, or  {\em expansion factor}  of the substitution.
 Also, under the assumption of  primitivity of $M_\sigma$,  the  shift $(X_\sigma,S)$   has particularly useful  ergodic  properties, such as unique ergodicity or minimality. For more details, see \cite{Queffelec:2010}.

 Analogously,  substitutions can  also  be  defined    on tiles, and iterating the substitution produces  self-affine tilings. 
  Tilings
are   considered in any dimension  and   over general combinatorial structures  with group actions.
Two   main  natural models of tilings are   Wang tiles  coming from  computer science   (see Section \ref{sec:comput}), and shifts coming from   dynamical  systems.
 In the tiling setting, we now work in a  $d$-dimensional geometric space, here ${\mathbb R}^d$,  rather than in  a symbolic  setting, with letters and words.

\begin{definition}[Self-affine tiling] \label{def:SA} Let $\{T_i; 1 \leq i \leq m\}$ be a set of   tiles (usually assumed to be polygons)  in   ${\mathbb R}^d$, called {\em prototiles}.
A  {\em self-affine tiling} ${\mathcal T} = \{ x_i + T_i ; x_i \in \Lambda_i, i \leq m\}$  in  ${\mathbb R}^d$  is   a tiling  fixed by a linear   expanding map  $\phi : {\mathbb R}^d \rightarrow {\mathbb R}^d$ that  maps every  prototile onto a union of prototiles ({\em expanding }means that all its eigenvalues are greater than one in modulus).  
The prototiles  are translated  by vectors in the sets  $\Lambda_i$ in order to form the tiling ${\mathcal T}$.
 With the map $\phi$ comes a  substitution rule $\Sigma$ that explains how to  divide a  tile  inflated by $\phi$  into a union of translated prototiles
and   an incidence matrix  $M_\Sigma$ which  indicates for each prototile $M$ how $\phi(M)$ is a union of prototiles.  More precisely,
there exist  finite sets $ {\mathcal D}_{ij} \subset {\mathbb R}^d$  ($1\leq i,j \leq k$) such that for any $j$ ($1\leq j \leq k)$
$$\Sigma(T_j)=\{ u +T_i,  u \in {\mathcal D}_{ij}, \  i\leq k\} \mbox{ with }
\phi  T_j=  \cup _{i=1}^k  (T_i + {\mathcal D}_{ij}).$$
Here all sets on the right-hand side must have disjoint interiors; it is also possible for
some of the sets ${\mathcal D}_{ij}$ to be empty. The matrix $M_\Sigma$ has as entries the cardinalities of the sets $ {\mathcal D}_{ij}$.

If  $M_\Sigma $  admits a strictly dominant eigenvalue   that is  positive, then as with substitutions, 
 this  eigenvalue  is called   the  inflation factor of the tiling. It equals the absolute value of the determinant of $\phi$, which controls the volume expansion.
The tiling  is said {\em self-similar} if $\phi$  is a similarity. 
\end{definition}
\begin{example}\label{ex:Fibo3}  We revisit  Example \ref{ex:Fibo} and \ref{ex:Fibo2}.  Here $d=1$. 
 We define a tiling with the Fibonacci substitution $\sigma: A\mapsto AB, B\mapsto A$ by starting with two intervals (prototiles), labelled respectively  $A$ and $B$, with lengths $\ell_A$ and $\ell_B$.  
To obtain a geometric realisation $\Sigma$ of the Fibonacci substitution $\sigma$ using two one-dimensional interval  tiles, we  thus  need to find tiles
 $A$ and $B$, with lengths $\ell_A$ and $\ell_B$, to be determined, where \begin{itemize}
 \item
 $\Sigma$ maps the tile $A$ to the two concatenated tiles  $AB$, 
 \item 
 $\Sigma$ maps the tile $B$ to the tile  $A$, and
 \item
 $\Sigma$ corresponds to an expansive map on $\mathbb R$.
 \end{itemize}
 The substitution $\Sigma$ maps $A$, whose length is $\ell_A$, to two tiles of combined length $\ell_A+\ell_B$,
 and the tile $B$ of length $\ell_B$ to a tile $A$ of length $\ell_A$. Since we want the substitution rule $\Sigma$ to act as a map that expands the  space with an inflation factor $\lambda$, this implies that
 \begin{align*}
 \lambda \ell_A = \ell_A+\ell_B \mbox{ and }  \lambda \ell_B = \ell_A ,
 \end{align*}
 i.e.,
  \begin{align*}
 \lambda (\ell_A ,\ell_B)= (\ell_A,\ell_B) \begin{pmatrix} 1&1\\1&0 \end{pmatrix} = (\ell_A,\ell_B) M_\Sigma ,
 \end{align*}
so that $\lambda$ must be a left eigenvector for $M_\Sigma$, and since $\varphi$ is the only expanding eigenvalue of $M_\Sigma$, $\lambda=\varphi$,  and we can take $\ell_A=\varphi$ and $\ell_B=1$.
We thus set  $A=[0, \varphi]$ and $B= [0, 1]$.  This defines a self-affine tiling (and even a self-similar one). The matrix $M_\Sigma$  equals the  matrix $M_{\sigma}$, and one has ${\mathcal D}_{11}= \{0\}$, ${\mathcal D}_{12}= \{\varphi\}$,
${\mathcal D}_{21}= \{0\}$, and ${\mathcal D}_{22}= \emptyset$.
   Iterating this tiling substitution rule gives us a geometric tiling of the line; see Figure \ref{Uwe}.
\end{example}

 Self-affine tilings   occur in   a wide  range  of  problems in dynamics; for instance, 
they are related to Markov partitions for hyperbolic maps and  radix representations  \cite{BerNum}. 
More generally, for tilings, not necessarily generated by substitutions, the associated  shift  dynamics is    given by  ${\mathbb R}^d$-translations 
yielding as   dynamical systems   orbit closures of  tilings  \cite{Rob04,Sol97}.
The closure is taken in the local topology   where  coincidences on a large ball around the origin up to a small translation  are considered.
 A  remarkable topological  property here  is that  the local structure of  the    tiling  space (i.e., the associated  dynamical system  endowed  with  the ${\mathbb R}^d$-action)    is   described as  the product of a Cantor set and a Euclidean space. Locally, it looks like  a disk crossed
 with a totally disconnected set.
 Tiling spaces  then  yield  particularly interesting  types of topological spaces   obtained as inverse limits of branched manifolds  with totally disconnected compact fibres \cite{Sadun}. This makes 
  cohomological methods  particularly  adapted in this setting,  and particularly the \u{C}ech cohomology  which gives information on how a  tiling can be deformed.      
This goes  beyond  the case of substitutive tilings. Consider  as  an illustration the works of A. Forrest, J. Hunton, and J. Kellendonk 
in  \cite{FHK:02bis}   for  cut and project  tilings  which, along the way,  illustrate   the efficiency of     topological methods      combined with  tools  from   $C^*$-algebras (see  also \cite{Putnam} for an elementary exposition). More generally,  topological  dynamics  and  topological notions   (e.g., maximal equicontinuous factor or  proximality relation) have led to beautiful developments, for instance   in  the computation of  the Ellis semigroup
\cite{ABKL:15,KY:20}.

\subsection{More on self-affine   tilings} \label{subsec:selfsim}
In this section, we revisit  Y. Meyer's Theorem  \ref{thm:meyer} in the setting  of self-similar tilings. According to Definition \ref{def:SA},  we consider   an expanding linear map $\phi : {\mathbb R}^d \rightarrow {\mathbb R}^d$ and  a  self-affine tiling   ${\mathcal T}$  fixed  by $\phi$.
Dynamicists ask which numbers  can appear as inflation factors, or which maps  can appear as  expanding maps,  of self-affine tilings, and which of these are inflation factors for substitutions that lead to quasicrystals. This has  led to  a vast literature, and 
 as in Theorem  \ref{thm:meyer},  it is articulated  in two steps:  firstly,  what are the   conditions   that the inflation factors should  satisfy? Secondly, can it be possible to realise  a   self-similar structure, and even a Meyer or a  model set,
with each  dilation  factor  satisfying  these  conditions?  We turn  to the  last question in Section  \ref{subsec:RF}.

We first  introduce two  arithmetic definitions.

 \begin{definition}[Perron number] 
 The algebraic integer $\beta$ is said to be  a 
    {\em Perron number}  if $|\beta|>|\alpha|$ for any other   of its  algebraic conjugates.    \end{definition}
  For example, the Perron--Frobenius theorem tells us that   the  dominant eigenvalue  of  a primitive non-negative matrix  is a Perron number.
   
    \begin{definition}[Pisot family] \label{def:PF}
Let  $\Lambda$  be  a  finite set of algebraic integers of modulus  larger than or equal to $1$. The set  $\Lambda$  is said to form a {\em Pisot family} if 
 for every $\lambda \in \Lambda$ and  for every  algebraic conjugate $ \lambda'$ of $\lambda$ with $|\lambda'| \geq 1$, then 
$\lambda' \in \Lambda$.

\end{definition}
 In  the one-dimensional case,   it follows from the work of D. Lind   \cite{Lind} that   $\lambda$ is an inflation factor of a self-similar  tiling if and only if  $\lambda $
  is a Perron number. 
  In the two-dimensional  self-similar  case (the dilation is  the same in every direction),   the   expansion map $\phi$, when viewed as acting on $\mathbb C$, equals $\lambda z$ where $\lambda$   is   a complex Perron number,
  as shown by  W. Thurston \cite{Thurston}.
   In the case where the   expansion map $\phi$ is diagonalizable,  R. Kenyon and B. Solomyak   \cite{KeSo:10}    \cite{Kenyon:PhD}
 showed   that   the  eigenvalues  of  $\phi$  are  algebraic integers. The general result  (without the diagonalizabiliy condition) was then established  by J.~Kwapisz
 \cite{Kwapisz:16} who proved that, with a suitable notion of Perron,  the map $\phi$ of a  self-affine tiling  is  integral algebraic and Perron.

 Now  that it was  established  what Perron numbers have to do with  being self-affine (or  even  being self-similar), Pisot  numbers entered into  play  when  further properties were required, such as some   of those shared by Meyer sets,   extending Y. Meyer's Theorem  \ref{thm:meyer}.  In this setting, the connection with  the notion of the  Pisot family    was     successfully developed by J.-Y. Lee and B. Solomyak  \cite{LeeSo:12}  for spectral considerations  in the case where the map $\phi$  is diagonalizable over the complex numbers, with all eigenvalues being algebraic conjugates of the same multiplicity.
They proved that the associated tiling shift has  relatively dense discrete spectrum  (i.e., the set of eigenvalues  has full rank)
 if and only if the spectrum of   the linear map $\phi$  forms a  Pisot family,   which  is  also equivalent to the   discrete set of control points for the tiling  to be  a Meyer set.

More generally,  we will see in Section \ref{subsec:Pisot}  how  Pisot   numbers   enter the picture  in this substitutive setting   not only  as inflation  factors, but  also   for model sets via  spectral  considerations.

\section{Pisot substitutions and aperiodic order} \label{sec:Pisot-Sturmian-Rauzy-Sadic}
     Now that inflation factors  are well  characterised, the question is to be able to  construct  self-similar tilings or point sets having this inflation factor  that are  also     Meyer  or even model sets.
There exist several   strategies for  such  realisations.  Among them, Rauzy fractals, discussed in  Section \ref{subsec:RF},   provide suitable windows for  explicit cut and project schemes.
But   before  discussing  them, we introduce  the  so-called  Pisot substitution conjecture.

\subsection{Spectral  properties  and the Pisot conjecture}  \label{subsec:Pisot}
 The spectral study of substitutive dynamical systems is a fundamental question (see for instance the monographs
  \cite{Queffelec:2010,PytheasFogg:2002}). This is   particularly  relevant  in the context of aperiodic order,
 as stressed by the     contribution  by  R. V. Moody in this volume  \cite{Moody:20}. Weak mixing, namely 
the absence of  eigenvalues, i.e., the absence of  Bragg peaks,  indicates a certain level of disorder   and conversely,   the existence of spectral eigenvalues
provides   dynamical factors which consist of group translations.     For instance,  constant length substitutions  yield  $p$-adic   factors.   In particular, 
discrete spectrum  deals with the possibility of providing substitutive  dynamical systems  with a  representation as a group translation, that is, a  dynamical system acting on a space  of  a geometric nature  (such as a torus for instance),  that is (measurably) isomorphic to the substitutive system under study.

 We have seen through the work of Y. Meyer   that the  connection with Pisot numbers arose     from the inception of  quasiorder, and this  has led to the following natural  definition.
 \begin{definition}[Pisot substitution]
 If  $\sigma$ is a word
substitution whose incidence matrix has a characteristic polynomial which is the  irreducible minimal polynomial of  a Pisot number,  then $\sigma$ is called a  {\em Pisot substitution.} 
\end{definition}
From    the eighties onwards,   Pisot  substitutions  have   attracted  much attention in the context of mathematical  quasicrystals.  We recall that the Fibonacci substitution provided   one of the   first  examples of a one-dimensional quasicrystal
  and more generally substitutions   were considered as promising examples of   one-dimensional quasicrystals; as we have seen, they 
     are   very  simple  algorithmic  rules   that   create configurations that display long-range order.   However  not  all substitutions  yield quasicrystals.
   E.  Bombieri  and J. E. Taylor  had  already asked  in \cite{BomTay:86,BomTay:87}   which substitutions 
 produce quasicrystals  and  they highlighted  the Pisot    algebraic restriction for  substitutions  to  yield pure discrete spectrum.   
 See also \cite{Sol97}   where  a self-similar tiling is proved  to admit a  discrete component
 in its spectrum if  and only if  its  inflation factor is a Pisot number.

  This culminated with the {\em Pisot  substitution  conjecture}, stated below for   substitutions defined on symbols. There exist various   formulations,  see the survey \cite{Akiyama&Barat&Berthe&Siegel:2004}.
 \begin{conjecture}[Pisot substitution  conjecture] \label{conj:Pisot} 
 Let $\sigma$  be   a Pisot  irreducible  substitution (i.e., the characteristic polynomial  of  its incidence matrix is
 the minimal polynomial of a Pisot number).  Then, 
 the shift $(X_{\sigma},S)$  has pure discrete spectrum, i.e., it is  isomorphic, in a  measure-theoretic sense,  to  
 a   translation on a compact abelian group.    
\end{conjecture}
In a nutshell,     the Pisot  arithmetic condition    induces  order, where order    is   expressed here in spectral and dynamical  terms
as  being  isomorphic to    the simplest  dynamical  systems, namely   group translations. 
The conjecture is  known to be true for the case  of substitutions on two letters   \cite{BD02,HS03}.     The difficulty for a larger alphabet  comes from  
 the  arithmetic of  higher degree  algebraic numbers,  which is more  difficult to manage than for quadratic  numbers.   Also, in  two dimensions, it is easier to identify the 
 occurrences of  so-called
 {\em coincidences}, the existence of which allows us to project  injectively  an infinite word  onto a group, here  the circle;  to do so, words are embedded   as  discrete  lines   in the plane and
coincidences  between  these  lines are to be found,  which is   easier in a plane (see Example \ref{ex:Tribo} for an illustration);
 generic behaviour is that once there is a coincidence, there are  infinitely many occurrences of the coincidence, and this implies that this projection is an injection.

\begin{theorem} \cite{BD02,HS03} \label{theo:2l}
Two-letter  Pisot substitutions have pure discrete spectrum.
\end{theorem}

The still open Pisot substitution conjecture, even if   solved in the  closely related  context of  beta-numeration  by M. Barge
  \cite{BargeBeta},
   shows that important parts of the picture are still  to be developed.
   Once again, one  particularly appealing  feature  concerning the works  developed   around   the Pisot substitution conjecture
   is that they involve   several  mathematical  approaches  and reformulations   involving  topological  dynamics, arithmetic,  combinatorics, fractal geometry,  etc.  See  as an illustration
\cite{Akiyama&Barat&Berthe&Siegel:2004} or \cite{thuswaldner2019boldsymbolsadic}.

     It follows from algebraic considerations, involving the  eigenvalues of the incidence  matrix, that Pisot substitution shifts and Pisot tiling dynamical systems must have a non-trivial rotation factor.
 Criteria for     obtaining   these dynamical eigenvalues  are   now well understood; see  for instance     \cite{Ferenczi-Mauduit-Nogueira,Sol:07}.   
 The difficult part  consists in  producing a measurable isomorphism. Fortunately,  there is a  wide range of  algorithmic   conditions, 
 called coincidence  conditions (such as briefly evoked above),
  for  checking  pure discrete spectrum, which date back to the  work of M. Dekking on constant-length substitutions, as described in   \cite{Queffelec:2010}.
   They occur for infinite words,  for higher-dimensional and for non-lattice based self-affine tilings,   and  even in the non-unimodular case. 
 See  for instance   \cite{Sol97,Lee&Moody:2001,SieThus:09,AkiLee:11}.   Also relevant is the notion of {\em almost-automorphy}, introduced by W. Veech \cite{Veech-1965}, which implies discrete spectrum for substitutions, and which incidentally was shown to imply spectral synthesis \cite{Veech-1969}, bringing us back again to Y. Meyer's work.

\subsection{Sturmian words and beyond   substitutions }\label{subsec:sturm}
As noted, the Pisot  substitution conjecture   has been proved for two--letter  irreducible Pisot substitutions (Theorem \ref{theo:2l}).    Such substitutions produce 
model sets.   Amongst them are  a  remarkable and   widely studied  family, that of the  Sturmian substitutions, which generate  shifts that belong to 
 the class of Sturmian shifts that we  now describe.
 
 Sturmian shifts are symbolic representations of irrational circle rotations. More  precisely,  consider the  rotation $R_{\alpha}$  acting  on $\mathbb{T}=\mathbb{R}/\mathbb{Z}$, with
 $R_{\alpha}(x)=x+\alpha \mod 1$, for $\alpha$  irrational.
 If we code the orbit of a point under the action of $R_{\alpha}$
using a   two-interval   partition of semi-open intervals  whose  lengths are respectively $\alpha$ and $1-\alpha$,   then  we obtain a {\em Sturmian word}, generating a {\em Sturmian shift} 
$X_{\alpha}$  consisting of 
 all  Sturmian words  of {\em angle} $\alpha$, that is,  all  Sturmian words coding $R_{\alpha}$.
  This    seminal  class  
of  symbolic dynamical systems,  introduced by M. Morse and  G. Hedlund   in \cite{Morse2},
   laid some of  the  foundations for  symbolic dynamics. For a thorough  description 
 of Sturmian  words,  see 
\cite{Berstel&Seebold:2002} and  \cite[Chapter 6]{PytheasFogg:2002}.

 The Fibonacci substitution   from Section \ref{subsec:subsdef} is the most cited example of a Sturmian substitution. When it is used to generate a tiling, then as we saw in Figure \ref{Uwe}
     the associated   point set is a  model set.
Sturmian words  have     provided  the  simplest examples of cut  and project sets:  
  the window  is here an interval    of unit  length and the associated tilings are  one-dimensional, with the cut and project scheme relying  on two lines
  as described  in Section \ref{sec:Pisot}.
Phrased in  a  more  geometric way,  they   code  discrete lines in  digital geometry  \cite{Rosenfeld}.
There is an impressive    literature devoted to the  study   of Sturmian words
and to  possible generalisations; let us mention for instance episturmian words, also called Arnoux-Rauzy words,
  which have  attracted a lot of attention
  \cite{Cassaigne-Ferenczi-Messaoudi:08,AHS:16bis}.
    
One key feature is  the  scale invariance  of Sturmian shifts. 
Not all Sturmian words are substitutive. Indeed the angle $\alpha$  of a substitutive  Sturmian shift $X_{\alpha}$   is  a quadratic number  \cite{Crisp} (see also  \cite{BEIR} for a  characterisation   of Sturmian words that are fixed points of  substitutions obtained in terms of Rauzy fractals discussed  in Section \ref{subsec:RF}).
 However  Sturmian words can all be generated  in terms of sequences of substitutions which we elaborate below.
  More precisely, Sturmian words   are   perfectly understood  and described  via 
   a representation  based on  substitutions and  Euclid's algorithm:
the  continued fraction expansion of    the angle  $\alpha$   provides   an infinite  product of  square matrices of size two with non-negative integer entries,
each of these matrices can be seen as  the  incidence  matrix of a  substitution, and   the action of a substitution  can   be seen
as a  combinatorial interpretation  of 
a  step of  Euclid's algorithm, as described below.
\begin{theorem} \cite{AR91}\label{theo:AR}
We consider  the substitutions  over the alphabet $\{0,1\}$ 
 defined by  $\sigma_0 \colon 0 \mapsto  01, 1 \mapsto 1$,  and $\sigma_1 \colon 0 \mapsto 0, 1 \mapsto 10$.
If $\alpha$  has  continued fraction expansion  $\alpha= [0; a_1+1, a_2, \cdots]$,  then  the Sturmian shift  $X_{\alpha}$ is generated by  the infinite word
$$ \lim_{n \rightarrow \infty} \sigma_0^{a_1} \sigma_1^{a_2}\sigma_0^{a_3} \sigma_1^{a_4} \cdots   \sigma_0^{a_{2n-1}} \sigma_1^{a_{2n}}  (0).$$
\end{theorem}

Seen through this lens, Sturmian words  are 
   described  using   a renormalization scheme  governed by continued fractions  via  the geodesic flow acting  on  the modular surface.
 This      crystallises     with   the study  of    interval exchanges in  relation with 
 the Teichm\"uller flow  through   the work, among others, of W. Veech,    H. Masur,  J.-C. Yoccoz, and  A. Avila (see e.g.    \cite{yoccoz,Bufetov:14}).

Theorem \ref{theo:AR} states  that  a Sturmian  word   is    the limit of  an infinite   composition of  these  substitutions;  such an approach  has    been formalised 
using the language of 
 {\em $S$-adic words},  which are infinite words generated when    a sequence $S$ of substitutions is applied, as opposed to when a single substitution is iterated   \cite{BerDel}.
The  $S$-adic  setting   pertains   to 
non-stationary   dynamics    (i.e., time inhomogeneous dynamics), which    consists in   working with  iterated  sequences of 
 transformations  drawn  according to a further dynamical system,  similarly as  for  random dynamics,     random  Markov chains and random  products of matrices.  This formalism allows us  to extend the  Pisot   conjecture  beyond algebraicity   \cite{BST:21}.
  The Pisot condition is then  replaced by the requirement 
    that the second Lyapunov exponent of the  infinite  associated  products of matrices  is negative.  
This  extended Pisot  conjecture    has been proved to hold for large relevant families of systems based on continued fractions expansions, such as the  Brun algorithm.
As  a striking outcome,   this yields  symbolic codings for almost every   translation of ${\mathbb T}^2$ \cite{BST:21}, paving the way for  the development of   equidistribution  results for
the associated two-dimensional  Kronecker sequences.   

\subsection{Rauzy fractals}   \label{subsec:RF} 

 In this section we describe a  construction  for word  substitutions  having a  Pisot  inflation factor  that generates  suitable and effective   windows  for  model sets;
this  generalises the Sturmian case for which  the window   in the internal  space is an interval. 
In fact, it is now well understood that a  substitution tiling with pure discrete spectrum  is   a  model set, such as described  by Lee in  \cite{Lee:07}.   She shows the equivalence
between   specific model sets, the so-called inter-model sets\footnote{Inter-model sets are model sets satisfying a topological condition that  is less restrictive than the
condition for the boundary of  having   zero measure   satisfied by  a regular model set.}, and pure point dynamical spectrum in the context of primitive substitution point sets,    where  the construction
 of the cut and project scheme involves an abstract internal space.  The  aim of this  section
is to   present   more explicit constructions of these internal spaces as
{\em Rauzy fractals}. This approach appeals in its combination of arithmetic and dynamics.

Recall  that for  $(X_{\sigma},S)$ to  have pure  discrete spectrum, it must be measurably  isomorphic to a group  translation.
 Thus, given    a Pisot substitution $\sigma$,  one wants     to provide a geometric representation for $(X_{\sigma},S)$
  as a translation on the torus, or more generally on a locally compact abelian group. And as a candidate for a fundamental domain for the expected translation, one associates with the substitutive dynamical system  $(X_{\sigma},S)$
  a   Rauzy fractal.

 Rauzy fractals were  first  introduced
in \cite{Rau82} in the case of   the so-called Tribonacci substitution in order to prove the following statement. 
\begin{theorem}\cite{Rau82} \label{theo:Tribo}
 Let $\sigma$   be  the   Tribonacci substitution defined over the alphabet ${\mathcal A}= \{A,B,C\}$  as as  $\sigma(A)=AB$, $\sigma(B)=AC$ and  $\sigma(C)=A$. 
The symbolic dynamical system $(X_{\sigma},S)$  is measure-theoretically  isomorphic
to the  translation $R_{\beta}$
 on the two-dimensional torus
${\mathbb T} ^2$  defined as 
$R_{\beta}:{\mathbb T}^2 \rightarrow {\mathbb T} ^2, \ 
x \mapsto x+(1/\beta,1/\beta ^2),$
where $\beta$, the Perron eigenvalue for $\sigma$, satisfies $\beta^3=\beta^2+\beta+1$.
\end{theorem}

\begin{figure}\label{Rauzy}
	\includegraphics[width=.4\textwidth]{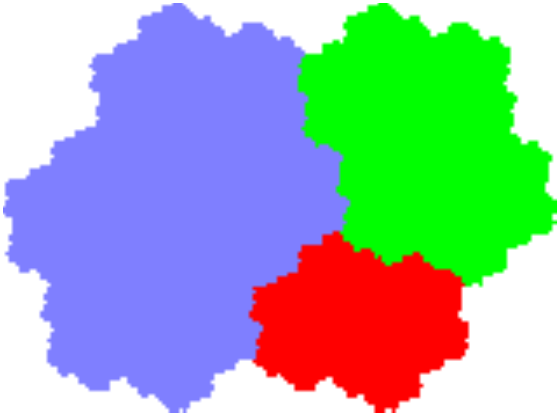}
	\caption{ The Rauzy fractal  }
\end{figure}

We now sketch how  to build  a fractal from the shift $(X_\sigma, S)$ as  a fundamental domain for the 
translation of Theorem \ref{theo:Tribo}.
\begin{example} \label{ex:Tribo}
Let $\alpha$ be one of the complex roots of $X^3-X^2-X-1$; we have $|\alpha|<1$. We will associate to a point $x\in X_\sigma$ an $\alpha$-expansion as follows.
The point  $x\in X_\sigma$ can  indeed be uniquely desubstituted, either as $x=\sigma(y)$, or as $x=S\sigma(y)$ for some $y\in X_\sigma$.    We can associate to $x$ an initial digit $\epsilon_0\in \{0,1\}$, depending on whether we shift or not. We repeat this procedure, desubstituting $y$ to obtain a second digit $\epsilon_1$. Iterating we obtain a sequence $(\epsilon_n)\in \{ 0,1\}^\mathbb N$, and it can be verified that  this sequence satisfies $\epsilon_n \epsilon_{n+1}\epsilon_{n+2}=0$ for each $n$.  We then project $x$ to  $z:=\sum_n \epsilon_n \alpha^n$;
 note the connection to beta-numeration  (evoked in Section \ref{numeration}) that appears; this connection 
  is studied in \cite{Thurston}.  Geometrically, this amounts   to the projection $\piint$   along the expanding eigendirection  of the matrix  $M_{\sigma}$ 
   with respect  to  the decomposition of  ${\mathbb R}^3$
  into   the  expanding eigendirection and  the contracting eigenplane of    the matrix  $M_{\sigma}$. 
The   image   $\piint(X_\sigma)$  of $X_\sigma$ is
the {\em Rauzy fractal $\mathcal R$}; see Figure \ref{Rauzy}. Depending on what we see at the 0-index $x_0$ of $x$,   we can further specify in which of the three regions $\mathcal R_0=\alpha \mathcal R$, $\mathcal R_1= \alpha^3+\alpha^2 \mathcal R$ and $\mathcal R_2= \alpha^3+\alpha^4 +\alpha^3 \mathcal R$ the point 
 $\piint(x)$ lives. The projection is injective outside of the boundaries between these three regions, and shifting inside $X_\sigma$ corresponds to an exchange from one region to another; see \cite{Mes98} for   details. The domain $\mathcal R$ is the internal space of a cut and project scheme which generates the Tribonacci tiling. \end{example}

Rauzy fractals can more generally   be associated with Pisot substitutions (see
\cite{AI01,BK06,Sie04}
 and the surveys
\cite{BS05b,PytheasFogg:2002}), as well as with Pisot beta-transformations and 
beta-shifts under the name of  central tiles \cite{Aki02}.
Rauzy fractals  can be  defined   in  a unimodular  case  \cite{BK06} (when  the  inflation   factor $\beta$  of a substitution  is a  Pisot unit)  as well as     in a $p$-adic setting \cite{Sie04,Sing:06}   (the primes $p$ that occur are  prime divisors of the norm of $\beta$ and   the  Rauzy fractal lives in 
 a finite product of  Euclidean  and $p$-adic spaces).   Moreover, the geometric properties of   Rauzy fractals  reflect the self-similarity of the associated substitutive dynamical system: they are solutions of (graph-directed) iterated function systems.  Note also  that  Rauzy fractals  have also  been studied as    quasicrystals
\cite{VM:01}.

Rauzy fractals do not only produce geometric representations of substitutive dynamical
  systems, but they also have  very interesting    Diophantine  applications. We refer to  
  \cite{Ito&Fujii&Higashino:2003,Adamczeswki&Frougny&Siegel&Steiner:2009}  for representative  examples.
Rauzy fractals    are  particularly  convenient for   providing   arithmetic descriptions
of  periodic orbits,
yielding   relevant generalisations of    Lagrange's and   Galois' theorems for  continued fractions.
 We  recall that  Lagrange's theorem states that
  the continued   fraction expansion  of $x $ is  eventually periodic
if and only if  either $x \in {\mathbb Q}$, or  $x$ is a  quadratic number. Furthermore, if 
 $x>1$ is  a quadratic  number, and $x'$  is its algebraic conjugate, then 
 the continued   fraction expansion  of  $x$ 
is   purely  periodic
if and only if  $x $   is    irrational  and  $ -1<x'<0$: this is   Galois' theorem.
Galois'  theorem can be proved dynamically by using the Gauss map  $T$ acting on $[0,1]$ as 
$x \mapsto \{1/x\}$.
A key idea used here   to provide a  dynamical proof of  Galois' theorem (and  thus  to describe purely periodic orbits)   is to    transform the map  $T$  that is  not one-to-one into an injective map  thanks  to a	 suitable realisation of its natural extension.

The same type of methods  can also  be applied   in   the  beta-numeration case allowing a  characterisation of  periodic orbits for a   Pisot number  $\beta$.
The condition  on  algebraic conjugates is then  expressed in terms of  Rauzy fractals, such as stated below, noting that Galois' theorem  expresses the pure periodicity of the continued fraction expansion of a quadratic   real number  $x$ in terms of  the location of  $x$ and its algebraic conjugate on the line. Here the location is expressed in terms of belonging to a Rauzy fractal.
\begin{theorem}
If  $\beta$ is a  Pisot number, then  $x\in[0,1]$   has an 
  eventually periodic  $\beta$-expansion  if and only if $x \in {\mathbb Q}(\beta)\cap[0,1]$.
Moreover, a   real number $ x \in {\mathbb Q}(\beta)\cap [0,1)$ has a purely periodic $\beta$-expansion
if and only if the vector composed of  $x$ and its conjugates belong to   the Rauzy  fractal  associated with $\beta$. \end{theorem}

The analogue of Lagrange's theorem has been  proved  in   \cite{BM77,Schmidt:80} and the  analogue  of Galois' theorem 
 in \cite{Akiyama&Barat&Berthe&Siegel:2004}
 (see also the references therein).
The proof  is here again based on an explicit realisation  of the natural  extension of the beta-transformation $T_{\beta}$.

\section{How  computability came to model sets} \label{sec:comput}

We have  seen  with Y. Meyer's Theorem \ref{thm:meyer}
 that algebraicity  enters naturally into play  for describing  self-similarity for   aperiodic order.
We will see in this section that  computability     is an increasingly   pertinent viewpoint  to consider.  The tiling  setting is  the one which is best suited  for this discussion 
with  the notion of a local rule,  which describes     the constraints on how adjacent  tiles can   be assembled.  More generally, local rules naturally model the energetic finite-range interactions between atoms.
Aperiodic order emerges when  local properties   enforce global order together with aperiodicity, and for tilings, this   results   in the research of  local  rules.

The notion of local rule   is  natural  when considering  Wang tiles. 
A  Wang  tile is a unit  square with  edges marked with symbols (or colors). 
A tiling by Wang tiles consists in assembling copies of the   tiles  (by translations only) so that symbols on shared   edges match.  These   tiles   were introduced by Wang   in 1961    for the study  of fragments of first order logics in the context of  the  {\em domino problem}, which     asks for the existence of an algorithm deciding whether a finite set of Wang tiles can tile the plane. 
The domino problem  is  known to be undecidable,  as proved by R. Berger in 1966.   In dynamical terms, this   amounts to ask whether
a given   two-dimensional shift  of finite type is  non-empty. This result  is based on   two ingredients:       firstly,  the simulation of Turing computations by tilings of the plane (usually via  local rules or substitutions), and      secondly, the existence of aperiodic  sets of tiles (in the sense that they   can only produce aperiodic tilings).  
Since the first examples of aperiodic shifts  of finite type  were  based on hierarchical structures~\cite{Berger:66,Robinson:71},
substitutive structures  have been  known  to  be able to  force  aperiodicity.
But more  than that,  in dimension $d \geq 2$, under natural assumptions, it is  known that
most  substitution tilings admit local rules (possibly after {\em decoration})  \cite{Mozes,Goodman:98}, as expressed in the following  statement  that  covers all known examples of hierarchical aperiodic tilings.

\begin{theorem}\cite{Goodman:98}
Every substitution tiling of  the $d$-dimensional Euclidean plane,  $d>1$, can be enforced with finite matching rules, subject to a mild condition. \end{theorem}

After  the connection between   higher-dimensional substitutions and local rules was established,
 the following question   became  natural: can we describe a tiling, which is not necessarily substitutive,  in  terms of   local rules? This question   emerged in mathematics as early as the 1990s  in the study of quasicrystals \cite{Katz,Levitov,Le:95,Socolar} (see more references in \cite{BedFer:15,FerSab:19}).
 This    has  known  particularly rich developments    for tilings   obtained by the  cut and project schemes.
 One important  feature here   is that it   comes down  both to computability   and algebraic considerations  on the associated parameters.

 Consider for instance 
 a    tiling   obtained  from a  $d$-dimensional space $E$   in ${\mathbb R}^n$, via the cut and project method.  When $d=1$ and $n=2$,
this yields a one-dimensional tiling, obtained by projecting a broken line  approximating   a real line embedded in a two-dimensional space.  This corresponds to  Sturmian words.  Note that they  cannot be defined 
 with local rules,     in  direct contrast  with  the higher-dimensional case, where the question is   how to force tiles to approximate a desired plane $E$ in ${\mathbb R}^n$  by specifying local rules.
The study of the  connections between the existence of local rules for  such  planar tilings and the parameters of its slope started with \cite{debruijn,Levitov,Le:95,Le:97,Socolar}.
  The first conditions   on  the parameters   of  cut and project  schemes  were  of an  algebraic nature.  In particular, it was proved in \cite{Le:97} that a slope enforced by undecorated local rules is necessarily algebraic (this is however not sufficient, see e.g. \cite{BedFer:15,BedFer:20}).  However, computability comes into play when the tiles can be decorated (i.e.,  in dynamical terms, when we  go from shifts of finite type  to    sofic ones). Decorations indeed allow  the  transfer  of information through the tiling, and this was used in \cite{FerSab:19} to prove  the following statement. \begin{theorem}\cite[Corollay 2]{FerSab:19}
A slope can be enforced by colored weak local rules if and only if its slope  is computable. 
\end{theorem} This  computational approach   is reminiscent of the proof of the undecidability of the domino problem. Colored local rules  are used  to encode simulations of Turing computations which check  that  only  planar  tilings  that  approximate  the  desired  planes  can  be  formed.  The fundamental tool is
   the use of so-called simulation theorems, which state that any effective one-dimensional  shift  can  be  obtained  as  the  subaction  of  a  two-dimensional  sofic shift  \cite{Hochman:09,AuSa:13,DRS}). 
The Penrose tiling is an example   where all these viewpoints  gather:  it can be considered simultaneously  as generated by a  substitution,   it is  a model set, thus   obtained via a cut and project scheme, and     it  can be endowed with local rules.

Nourished by this computability   viewpoint,  aperiodic order    has known particularly  fruitful  developments  with the study of tilings,   spreading also in the direction of higher-dimensional 
  word combinatorics, symbolic dynamics    or else  group combinatorics,   with   symbolic dynamics on  groups and the domino problem for them.
  In this  latter setting,   the decidability of the domino problem is reinterpreted as a group property,   inspired by Higman's embedding theorem 
which states  that any finitely generated and recursively presented group  can be embedded  as a subgroup of a finitely presented group,
 such as   developed in   \cite{Jean:16,JV:19}. 
 
  \bigskip
  
 In conclusion,  Y. Meyer's Theorem  \ref{thm:meyer} and its mathematical  successors    show that Pisot numbers, and more generally   arithmetic and Diophantine  approximation,
   are natural actors of  aperiodic order   when self-similarity occurs.   As an illustration, observe  that the most classical   combinatorial and dynamical  measures    of disorder    developed for   cut and project  structures  can be evaluated in terms of Diophantine properties   such as  
developed e.g.  in \cite{HKWS:16,HJKW:19}   with  in particular  nice connections with the Littlewood conjecture  \cite{HKW:18}.
However, the  last developments  of  the study of  quasiorder   show  that  computability is  also a   viewpoint   that has to be considered, with    the possibility   of encoding  some computation in tilings. We have seen that  this  has   emerged as early as in  the 1960s with the introduction of Wang tiles.
 This  concept  then crystallised in a striking way  50 years later with the work of M.~Hochman \cite{Hochman:09}     and  with the characterisation, obtained by M.~Hochman and  T.~Meyerovitch, of the entropies of multidimensional finite type shifts. If these entropies are  logarithms of Perron numbers in  
the one-dimensional case, M. Hochman and T.  Meyerovitch    made a fundamental change of perspective  by  proving that  in the higher-dimensional case, 
these entropies  are characterised as   the  recursive numbers
that are right  recursively  enumerable \cite{HochmanMeyer}.

\bibliographystyle{alpha}
\bibliography{Meyerbib}

\begin{thebibliography}{HKWS16}

\bibitem[ABBS08]{Akiyama&Barat&Berthe&Siegel:2004}
S.~Akiyama, G.~Barat, V.~Berth\'{e}, and A.~Siegel.
\newblock Boundary of central tiles associated with {P}isot beta-numeration and
  purely periodic expansions.
\newblock {\em Monatsh. Math.}, 155(3-4):377--419, 2008.

\bibitem[ABF05]{ABF:05}
S.~Akiyama, F.~Bassino, and C.~Frougny.
\newblock Arithmetic {M}eyer sets and finite automata.
\newblock {\em Inform. and Comput.}, 201(2):199--215, 2005.

\bibitem[ABKL15]{ABKL:15}
J.-B. Aujogue, M.~Barge, J.~Kellendonk, and D.~Lenz.
\newblock Equicontinuous factors, proximality and {E}llis semigroup for
  {D}elone sets.
\newblock In {\em Mathematics of aperiodic order}, volume 309 of {\em Progr.
  Math.}, pages 137--194. Birkh\"{a}user/Springer, Basel, 2015.

\bibitem[AFSS10]{Adamczeswki&Frougny&Siegel&Steiner:2009}
B.~Adamczewski, Ch. Frougny, A.~Siegel, and W.~Steiner.
\newblock Rational numbers with purely periodic beta-expansion.
\newblock {\em Bull. London Math. Soc.}, 42:538--552, 2010.

\bibitem[AHS16]{AHS:16bis}
A.~Avila, P.~Hubert, and A.~Skripchenko.
\newblock On the {H}ausdorff dimension of the {R}auzy gasket.
\newblock {\em Bull. Soc. Math. France}, 144:539--568, 2016.

\bibitem[AI01]{AI01}
P.~Arnoux and S.~Ito.
\newblock Pisot substitutions and {R}auzy fractals.
\newblock {\em Bull. Belg. Math. Soc. Simon Stevin}, 8(2):181--207, 2001.

\bibitem[AK13]{AK:13}
S.~Akiyama and V.~Komornik.
\newblock Discrete spectra and {P}isot numbers.
\newblock {\em J. Number Theory}, 133(2):375--390, 2013.

\bibitem[Aki02]{Aki02}
S.~Akiyama.
\newblock On the boundary of self affine tilings generated by {P}isot numbers.
\newblock {\em J. Math. Soc. Japan}, 54(2):283--308, 2002.

\bibitem[AL11]{AkiLee:11}
S.~Akiyama and J.-Y. Lee.
\newblock Algorithm for determining pure pointedness of self-affine tilings.
\newblock {\em Adv. Math.}, 226(4):2855--2883, 2011.

\bibitem[AR91]{AR91}
P.~Arnoux and G.~Rauzy.
\newblock Repr\'esentation g\'eom\'etrique de suites de complexit\'e $2n+1$.
\newblock {\em Bull. Soc. Math. France}, 119(2):199--215, 1991.

\bibitem[AS03]{ASBook}
J.-P. Allouche and J.~Shallit.
\newblock {\em Automatic sequences}.
\newblock Cambridge University Press, Cambridge, 2003.
\newblock Theory, applications, generalizations.

\bibitem[AS13]{AuSa:13}
N.~Aubrun and M.~Sablik.
\newblock Simulation of effective subshifts by two-dimensional subshifts of
  finite type.
\newblock {\em Acta Appl. Math.}, 126:35--63, 2013.

\bibitem[Bar18]{BargeBeta}
M.~Barge.
\newblock The {P}isot conjecture for {$\beta$}-substitutions.
\newblock {\em Ergodic Theory Dynam. Systems}, 38:444--472, 2018.

\bibitem[BBLT06]{BerNum}
G.~Barat, V.~Berth\'{e}, P.~Liardet, and J.~Thuswaldner.
\newblock Dynamical directions in numeration.
\newblock {\em Ann. Inst. Fourier (Grenoble)}, 56:1987--2092, 2006.
\newblock Num\'{e}ration, pavages, substitutions.

\bibitem[BD02]{BD02}
M.~Barge and B.~Diamond.
\newblock Coincidence for substitutions of {P}isot type.
\newblock {\em Bull. Soc. Math. France}, 130:619--626, 2002.

\bibitem[BD14]{BerDel}
V.~Berth\'{e} and V.~Delecroix.
\newblock Beyond substitutive dynamical systems: {$S$}-adic expansions.
\newblock In {\em Numeration and substitution 2012}, RIMS K\^{o}ky\^{u}roku
  Bessatsu, B46, pages 81--123. Res. Inst. Math. Sci. (RIMS), Kyoto, 2014.

\bibitem[BEIR07]{BEIR}
V.~Berth\'{e}, H.~Ei, S.~Ito, and H.~Rao.
\newblock On substitution invariant {S}turmian words: an application of {R}auzy
  fractals.
\newblock {\em Theor. Inform. Appl.}, 41(3):329--349, 2007.

\bibitem[Ber66]{Berger:66}
R.~Berger.
\newblock The undecidability of the domino problem.
\newblock {\em Mem. Amer. Math. Soc.}, 66:72, 1966.

\bibitem[Ber77]{BM77}
A.~Bertrand.
\newblock D\'{e}veloppements en base de {P}isot et r\'{e}partition modulo
  {$1$}.
\newblock {\em C. R. Acad. Sci. Paris S\'{e}r. A-B}, 285(6):A419--A421, 1977.

\bibitem[BF15]{BedFer:15}
N.~B\'{e}daride and Th. Fernique.
\newblock When periodicities enforce aperiodicity.
\newblock {\em Comm. Math. Phys.}, 335(3):1099--1120, 2015.

\bibitem[BF20]{BedFer:20}
N.~B\'{e}daride and Th. Fernique.
\newblock Canonical projection tilings defined by patterns.
\newblock {\em Geom. Dedicata}, 208:157--175, 2020.

\bibitem[BFGK98]{BFGK:98}
\v{C}. Burd\'{\i}k, Ch. Frougny, J.~P. Gazeau, and R.~Krejcar.
\newblock Beta-integers as natural counting systems for quasicrystals.
\newblock {\em J. Phys. A}, 31(30):6449--6472, 1998.

\bibitem[BG13]{AperBook}
M.~Baake and U.~Grimm.
\newblock {\em Aperiodic order. {V}ol. 1}, volume 149 of {\em Encyclopedia of
  Mathematics and its Applications}.
\newblock Cambridge University Press, Cambridge, 2013.

\bibitem[BK06]{BK06}
M.~Barge and J.~Kwapisz.
\newblock Geometric theory of unimodular {P}isot substitutions.
\newblock {\em Amer. J. Math.}, 128:1219--1282, 2006.

\bibitem[BM00]{CRM:2000}
M.~Baake and R.~V. Moody, editors.
\newblock {\em Directions in mathematical quasicrystals}, volume~13 of {\em CRM
  Monograph Series}.
\newblock American Mathematical Society, Providence, RI, 2000.

\bibitem[BS02]{Berstel&Seebold:2002}
J.~Berstel and P.~S\'e\'ebold.
\newblock Sturmian words.
\newblock In M.~Lothaire, editor, {\em Algebraic Combinatorics on Words},
  volume~90 of {\em Encyclopedia of Mathematics and Its Applications}, pages
  45--110. Cambridge University Press, 2002.

\bibitem[BS05]{BS05b}
V.~Berth{\'e} and A.~Siegel.
\newblock Tilings associated with beta-numeration and substitutions.
\newblock {\em Integers}, 5(3):A2, 46 pp. (electronic), 2005.

\bibitem[BST23]{BST:21}
V.~Berth\'{e}, W.~Steiner, and J.~M. Thuswaldner.
\newblock Multidimensional continued fractions and symbolic codings of toral
  translations.
\newblock {\em J. Eur. Math. Soc. (JEMS)}, 25(12):4997--5057, 2023.

\bibitem[BT86]{BomTay:86}
E.~Bombieri and J.~E. Taylor.
\newblock Which distributions of matter diffract? {A}n initial investigation.
\newblock volume~47, pages C3--19--C3--28. 1986.
\newblock International workshop on aperiodic crystals (Les Houches, 1986).

\bibitem[BT87]{BomTay:87}
E.~Bombieri and J.~E. Taylor.
\newblock Quasicrystals, tilings, and algebraic number theory: some preliminary
  connections.
\newblock In {\em The legacy of {S}onya {K}ovalevskaya ({C}ambridge, {M}ass.,
  and {A}mherst, {M}ass., 1985)}, volume~64 of {\em Contemp. Math.}, pages
  241--264. Amer. Math. Soc., Providence, RI, 1987.

\bibitem[Buf14]{Bufetov:14}
A.~I. Bufetov.
\newblock Limit theorems for translation flows.
\newblock {\em Ann. of Math. (2)}, 179:431--499, 2014.

\bibitem[Bug96]{Bugeaud:96}
Y.~Bugeaud.
\newblock On a property of {P}isot numbers and related questions.
\newblock {\em Acta Math. Hungar.}, 73(1-2):33--39, 1996.

\bibitem[CFM08]{Cassaigne-Ferenczi-Messaoudi:08}
J.~Cassaigne, S.~Ferenczi, and A.~Messaoudi.
\newblock Weak mixing and eigenvalues for {A}rnoux-{R}auzy sequences.
\newblock {\em Ann. Inst. Fourier (Grenoble)}, 58(6):1983--2005, 2008.

\bibitem[CMPS93]{Crisp}
D.~Crisp, W.~Moran, A.~Pollington, and P.~Shiue.
\newblock Substitution invariant cutting sequences.
\newblock {\em J. Th\'{e}or. Nombres Bordeaux}, 5(1):123--137, 1993.

\bibitem[dB81]{debruijn}
N.~G. de~Bruijn.
\newblock Algebraic theory of {P}enrose's nonperiodic tilings of the plane.
  {I}, {II}.
\newblock {\em Nederl. Akad. Wetensch. Indag. Math.}, 43(1):39--52, 53--66,
  1981.

\bibitem[DRS10]{DRS}
B.~Durand, A.~Romashchenko, and A.~Shen.
\newblock Effective closed subshifts in 1{D} can be implemented in 2{D}.
\newblock In {\em Fields of logic and computation}, volume 6300 of {\em Lecture
  Notes in Comput. Sci.}, pages 208--226. Springer, Berlin, 2010.

\bibitem[EJK90]{EJK:90}
P.~Erd\"{o}s, I.~Jo\'{o}, and V.~Komornik.
\newblock Characterization of the unique expansions
  {$1=\sum^\infty_{i=1}q^{-n_i}$} and related problems.
\newblock {\em Bull. Soc. Math. France}, 118(3):377--390, 1990.

\bibitem[EJK98]{EJK:98}
P.~Erd\"{o}s, I.~Jo\'{o}, and V.~Komornik.
\newblock On the sequence of numbers of the form
  {$\epsilon_0+\epsilon_1q+\cdots+\epsilon_nq^n,\ \epsilon_i\in\{0,1\}$}.
\newblock {\em Acta Arith.}, 83(3):201--210, 1998.

\bibitem[EK98]{EK}
P.~Erd\H{o}s and V.~Komornik.
\newblock Developments in non-integer bases.
\newblock {\em Acta Math. Hungar.}, 79(1-2):57--83, 1998.

\bibitem[FHK02]{FHK:02bis}
A.~Forrest, J.~Hunton, and J.~Kellendonk.
\newblock Topological invariants for projection method patterns.
\newblock {\em Mem. Amer. Math. Soc.}, 159(758), 2002.

\bibitem[FMN96]{Ferenczi-Mauduit-Nogueira}
S.~Ferenczi, C.~Mauduit, and A.~Nogueira.
\newblock Substitution dynamical systems: algebraic characterization of
  eigenvalues.
\newblock {\em Ann. Sci. \'{E}cole Norm. Sup. (4)}, 29(4):519--533, 1996.

\bibitem[FS92]{FrougnySolomyak92}
C.~Frougny and B.~Solomyak.
\newblock Finite beta-expansions.
\newblock {\em Ergodic Theory Dynamical Systems}, 12:45--82, 1992.

\bibitem[FS19]{FerSab:19}
Th. Fernique and M.~Sablik.
\newblock Weak colored local rules for planar tilings.
\newblock {\em Ergodic Theory Dynam. Systems}, 39(12):3322--3346, 2019.

\bibitem[FW02]{FW:02}
D-J. Feng and Z.-Y. Wen.
\newblock A property of {P}isot numbers.
\newblock {\em J. Number Theory}, 97(2):305--316, 2002.

\bibitem[Gar62]{Garsia}
A.~M. Garsia.
\newblock Arithmetic properties of {B}ernoulli convolutions.
\newblock {\em Trans. Amer. Math. Soc.}, 102:409--432, 1962.

\bibitem[GMP03]{GMP:03}
L.-S. Guimond, Z.~Mas\'{a}kov\'{a}, and E.~Pelantov\'{a}.
\newblock Combinatorial properties of infinite words associated with
  cut-and-project sequences.
\newblock {\em J. Th\'{e}or. Nombres Bordeaux}, 15(3):697--725, 2003.

\bibitem[GMP06]{GMP:06}
J.-P. Gazeau, Z.~Mas\'{a}kov\'{a}, and E.~Pelantov\'{a}.
\newblock Nested quasicrystalline discretisations of the line.
\newblock In {\em Physics and number theory}, volume~10 of {\em IRMA Lect.
  Math. Theor. Phys.}, pages 79--131. Eur. Math. Soc., Z\"{u}rich, 2006.

\bibitem[GQ]{GQ:20}
D.~Gratias and M.~Quiquandon.
\newblock Cristallographie, quasicristaux et {Y}ves {M}eyer.
\newblock This volume.

\bibitem[GS98]{Goodman:98}
C.~Goodman-Strauss.
\newblock Matching rules and substitution tilings.
\newblock {\em Ann. of Math. (2)}, 147(1):181--223, 1998.

\bibitem[Had98]{Hadamard}
J.~Hadamard.
\newblock Sur la forme des lignes g\'{e}od\'{e}siques \`a l'infini et sur les
  g\'{e}od\'{e}siques des surfaces r\'{e}gl\'{e}es du second ordre.
\newblock {\em Bull. Soc. Math. France}, 26:195--216, 1898.

\bibitem[HJKW19]{HJKW:19}
A.~Haynes, A.~Julien, H.~Koivusalo, and J.~Walton.
\newblock Statistics of patterns in typical cut and project sets.
\newblock {\em Ergodic Theory Dynam. Systems}, 39(12):3365--3387, 2019.

\bibitem[HKW18]{HKW:18}
A.~Haynes, H.~Koivusalo, and J.~Walton.
\newblock Perfectly ordered quasicrystals and the {L}ittlewood conjecture.
\newblock {\em Trans. Amer. Math. Soc.}, 370:4975--4992, 2018.

\bibitem[HKWS16]{HKWS:16}
A.~Haynes, H.~Koivusalo, J.~Walton, and L.~Sadun.
\newblock Gaps problems and frequencies of patches in cut and project sets.
\newblock {\em Math. Proc. Cambridge Philos. Soc.}, 161:65--85, 2016.

\bibitem[HM10]{HochmanMeyer}
M.~Hochman and T.~Meyerovitch.
\newblock A characterization of the entropies of multidimensional shifts of
  finite type.
\newblock {\em Ann. of Math. (2)}, 171:2011--2038, 2010.

\bibitem[Hoc09]{Hochman:09}
M.~Hochman.
\newblock On the dynamics and recursive properties of multidimensional symbolic
  systems.
\newblock {\em Invent. Math.}, 176:131--167, 2009.

\bibitem[HP19]{Abel}
H.~Holden and R.~Piene, editors.
\newblock {\em The {A}bel {P}rize 2013--2017}.
\newblock Springer, Cham, 2019.

\bibitem[HS03]{HS03}
M.~Hollander and B.~Solomyak.
\newblock Two-symbol {P}isot substitutions have pure discrete spectrum.
\newblock {\em Ergodic Theory Dynam. Systems}, 23:533--540, 2003.

\bibitem[IFHY03]{Ito&Fujii&Higashino:2003}
S.~Ito, J.~Fujii, H.~Higashino, and S.-i. Yasutomi.
\newblock On simultaneous approximation to {$(\alpha,\alpha^2)$} with
  {$\alpha^3+k\alpha-1=0$}.
\newblock {\em J. Number Theory}, 99(2):255--283, 2003.

\bibitem[Jea16]{Jean:16}
E.~Jeandel.
\newblock Computability in symbolic dynamics.
\newblock In {\em Pursuit of the universal}, volume 9709 of {\em Lecture Notes
  in Comput. Sci.}, pages 124--131. Springer, 2016.

\bibitem[JV19]{JV:19}
E.~Jeandel and P.~Vanier.
\newblock A characterization of subshifts with computable language.
\newblock In {\em 36th {I}nternational {S}ymposium on {T}heoretical {A}spects
  of {C}omputer {S}cience}, volume 126 of {\em LIPIcs. Leibniz Int. Proc.
  Inform.}, pages Art. No. 40, 16. Schloss Dagstuhl. Leibniz-Zent. Inform.,
  Wadern, 2019.

\bibitem[Kat95]{Katz}
A.~Katz.
\newblock Matching rules and quasiperiodicity: the octagonal tilings.
\newblock In {\em Beyond quasicrystals ({L}es {H}ouches, 1994)}, pages
  141--189. Springer, Berlin, 1995.

\bibitem[Ken90]{Kenyon:PhD}
R.~W. Kenyon.
\newblock {\em Self-similar tilings}.
\newblock ProQuest LLC, Ann Arbor, MI, 1990.
\newblock Thesis (Ph.D.)--Princeton University.

\bibitem[KLS15]{MAP:15}
J.~Kellendonk, D.~Lenz, and J.~Savinien, editors.
\newblock {\em Mathematics of aperiodic order}, volume 309 of {\em Progress in
  Mathematics}.
\newblock Birkh\"{a}user/Springer, Basel, 2015.

\bibitem[KS10]{KeSo:10}
R.~Kenyon and B.~Solomyak.
\newblock On the characterization of expansion maps for self-affine tilings.
\newblock {\em Discrete Comput. Geom.}, 43(3):577--593, 2010.

\bibitem[Kwa16]{Kwapisz:16}
J.~Kwapisz.
\newblock Inflations of self-affine tilings are integral algebraic {P}erron.
\newblock {\em Invent. Math.}, 205(1):173--220, 2016.

\bibitem[KY20]{KY:20}
J.~Kellendonk and R.~Yassawi.
\newblock The {E}llis semigroup of bijective substitutions.
\newblock {\em Groups, Geometry and Dynamics}, 2020.

\bibitem[Lag96]{Lagarias:96}
J.~C. Lagarias.
\newblock Meyer's concept of quasicrystal and quasiregular sets.
\newblock {\em Comm. Math. Phys.}, 179(2):365--376, 1996.

\bibitem[Lag99]{Lagarias:99}
J.~C. Lagarias.
\newblock Geometric models for quasicrystals {I}. {D}elone sets of finite type.
\newblock {\em Discrete Comput. Geom.}, 21(2):161--191, 1999.

\bibitem[Le95]{Le:95}
T.~Q.~T. Le.
\newblock Local rules for pentagonal quasi-crystals.
\newblock {\em Discrete Comput. Geom.}, 14(1):31--70, 1995.

\bibitem[Le97]{Le:97}
T.~Q.~T. Le.
\newblock Local rules for quasiperiodic tilings.
\newblock In {\em The mathematics of long-range aperiodic order ({W}aterloo,
  {ON}, 1995)}, volume 489 of {\em NATO Adv. Sci. Inst. Ser. C Math. Phys.
  Sci.}, pages 331--366. Kluwer Acad. Publ., Dordrecht, 1997.

\bibitem[Lee07]{Lee:07}
J.-Y. Lee.
\newblock Substitution {D}elone sets with pure point spectrum are inter-model
  sets.
\newblock {\em J. Geom. Phys.}, 57(11):2263--2285, 2007.

\bibitem[Lev88]{Levitov}
L.~S. Levitov.
\newblock Local rules for quasicrystals.
\newblock {\em Comm. Math. Phys.}, 119(4):627--666, 1988.

\bibitem[Lin84]{Lind}
D.~A. Lind.
\newblock The entropies of topological {M}arkov shifts and a related class of
  algebraic integers.
\newblock {\em Ergodic Theory Dynam. Systems}, 4(2):283--300, 1984.

\bibitem[LM01]{Lee&Moody:2001}
J.-Y. Lee and R.~V. Moody.
\newblock Lattice substitution systems and model sets.
\newblock {\em Discrete Comput. Geom.}, 25(2):173--201, 2001.

\bibitem[Lot02]{Lothaire:2002}
M.~Lothaire.
\newblock {\em Algebraic combinatorics on words}, volume~90 of {\em
  Encyclopedia of Mathematics and its Applications}.
\newblock Cambridge University Press, Cambridge, 2002.

\bibitem[LS12]{LeeSo:12}
J.-Y. Lee and B.~Solomyak.
\newblock Pisot family self-affine tilings, discrete spectrum, and the {M}eyer
  property.
\newblock {\em Discrete Contin. Dyn. Syst.}, 32(3):935--959, 2012.

\bibitem[Mes98]{Mes98}
A.~Messaoudi.
\newblock Propri\'et\'es arithm\'etiques et dynamiques du fractal de {R}auzy.
\newblock {\em J. Th\'eor. Nombres Bordeaux}, 10(1):135--162, 1998.

\bibitem[Mey70]{Meyer:70}
Y.~Meyer.
\newblock {\em Nombres de {P}isot, nombres de {S}alem et analyse harmonique}.
\newblock Lecture Notes in Mathematics, Vol. 117. Springer-Verlag, Berlin-New
  York, 1970.
\newblock Cours Peccot donn\'{e} au Coll\`ege de France en avril-mai 1969.

\bibitem[Mey72]{Meyer:72}
Y.~Meyer.
\newblock {\em Algebraic numbers and harmonic analysis}.
\newblock North-Holland Publishing Co., Amsterdam-London; American Elsevier
  Publishing Co., Inc., New York, 1972.
\newblock North-Holland Mathematical Library, Vol. 2.

\bibitem[Mey95]{Meyer:95}
Y.~Meyer.
\newblock Quasicrystals, {D}iophantine approximation and algebraic numbers.
\newblock In {\em Beyond quasicrystals ({L}es {H}ouches, 1994)}, pages 3--16.
  Springer, Berlin, 1995.

\bibitem[Mey12]{Meyer:12}
Y.~Meyer.
\newblock Quasicrystals, almost periodic patterns, mean-periodic functions and
  irregular sampling.
\newblock {\em Afr. Diaspora J. Math.}, 13(1):1--45, 2012.

\bibitem[Mey20]{Meyer:20}
Y.~Meyer.
\newblock From {S}alomon {B}ochner to {D}an {S}hechtman.
\newblock {\em Trans. R. Norw. Soc. Sci. Lett.}, 1:1--22, 2020.

\bibitem[MH40]{Morse2}
M.~Morse and G.~A. Hedlund.
\newblock Symbolic dynamics {I}{I}. {S}turmian trajectories.
\newblock {\em Amer. J. Math.}, 62:1--42, 1940.

\bibitem[Moo]{Moody:20}
R.~V. Moody.
\newblock Meyer sets and diffraction.
\newblock This volume.

\bibitem[Moo97]{Moody:97}
R.~V. Moody.
\newblock Meyer sets and their duals.
\newblock In {\em The mathematics of long-range aperiodic order ({W}aterloo,
  {ON}, 1995)}, volume 489 of {\em NATO Adv. Sci. Inst. Ser. C Math. Phys.
  Sci.}, pages 403--441. Kluwer Acad. Publ., Dordrecht, 1997.

\bibitem[Moo00]{Moody:00}
R.~V. Moody.
\newblock Model sets: a survey.
\newblock In {\em From Quasicrystals to More Complex Systems}. Springer Verlag,
  2000.
\newblock F. Axel, F. D\'enoyer, and J.-P. Gazeau, Centre de physique Les
  Houches.

\bibitem[Mor21]{Morse}
H.~M. Morse.
\newblock Recurrent geodesics on a surface of negative curvature.
\newblock {\em Trans. Amer. Math. Soc.}, 22, 1921.

\bibitem[Moz89]{Mozes}
S.~Mozes.
\newblock Tilings, substitution systems and dynamical systems generated by
  them.
\newblock {\em J. Analyse Math.}, 53:139--186, 1989.

\bibitem[OU]{OU:20}
A.~Olevskii and A.~Ulanovskii.
\newblock Meyer sets and related problems.
\newblock This volume.

\bibitem[Par60]{Parry}
W.~Parry.
\newblock On the {$\beta $}-expansions of real numbers.
\newblock {\em Acta Math. Acad. Sci. Hungar.}, 11:401--416, 1960.

\bibitem[PM07]{PelMas:07}
E.~Pelantov\'{a} and Z.~Mas\'{a}kov\'{a}.
\newblock Quasicrystals: algebraic, combinatorial and geometrical aspects.
\newblock In {\em Physics and theoretical computer science}, volume~7 of {\em
  NATO Secur. Sci. Ser. D Inf. Commun. Secur.}, pages 113--131. IOS, Amsterdam,
  2007.

\bibitem[Put18]{Putnam}
I.~F. Putnam.
\newblock {\em Cantor minimal systems}, volume~70 of {\em University Lecture
  Series}.
\newblock American Mathematical Society, Providence, RI, 2018.

\bibitem[{Pyt}02]{PytheasFogg:2002}
N.~{Pytheas Fogg}.
\newblock {\em Substitutions in Dynamics, Arithmetics and Combinatorics},
  volume 1794 of {\em Lecture Notes in Mathematics}.
\newblock Springer Verlag, 2002.
\newblock Ed. by V. {Berth\'e} and S. Ferenczi and C. Mauduit and A. Siegel.

\bibitem[Que10]{Queffelec:2010}
M.~Queff{\'e}lec.
\newblock {\em Substitution dynamical systems---spectral analysis}, volume 1294
  of {\em Lecture Notes in Mathematics}.
\newblock Springer-Verlag, Berlin, second edition, 2010.

\bibitem[Rau82]{Rau82}
G.~Rauzy.
\newblock Nombres alg\'ebriques et substitutions.
\newblock {\em Bull. Soc. Math. France}, 110(2):147--178, 1982.

\bibitem[R{\'e}n57]{Renyi}
A.~R{\'e}nyi.
\newblock Representations for real numbers and their ergodic properties.
\newblock {\em Acta Math. Acad. Sci. Hungar.}, 8:477--493, 1957.

\bibitem[RK01]{Rosenfeld}
A.~Rosenfeld and R.~Klette.
\newblock Digital straightness.
\newblock {\em Electronic Notes in Theoretical Computer Science}, 46:1 -- 32,
  2001.
\newblock IWCIA 2001, 8th International Workshop on Combinatorial Image
  Analysis.

\bibitem[Rob71]{Robinson:71}
R.~M. Robinson.
\newblock Undecidability and nonperiodicity for tilings of the plane.
\newblock {\em Invent. Math.}, 12:177--209, 1971.

\bibitem[Rob04]{Rob04}
E.~A. Robinson, Jr.
\newblock Symbolic dynamics and tilings of {${{\mathbb R}^d}$}.
\newblock In {\em Symbolic dynamics and its applications}, volume~60 of {\em
  Proc. Sympos. Appl. Math., Amer. Math. Soc. Providence, RI}, pages 81--119,
  2004.

\bibitem[Sad08]{Sadun}
L.~Sadun.
\newblock {\em Topology of tiling spaces}, volume~46 of {\em University Lecture
  Series}.
\newblock American Mathematical Society, Providence, RI, 2008.

\bibitem[SBGC84]{QC}
D.~Shechtman, I.~Blech, D.~Gratias, and J.~W. Cahn.
\newblock Metallic phase with long-range orientational order and no
  translational symmetry.
\newblock {\em Phys. Rev. Lett.}, 53:1951--1953, Nov 1984.

\bibitem[Sch80]{Schmidt:80}
K.~Schmidt.
\newblock On periodic expansions of {P}isot numbers and {S}alem numbers.
\newblock {\em Bull. London Math. Soc.}, 12(4):269--278, 1980.

\bibitem[Sen95]{sene}
M.~Senechal.
\newblock {\em Quasicrystals and geometry}.
\newblock Cambridge University Press, Cambridge, 1995.

\bibitem[Sie04]{Sie04}
A.~Siegel.
\newblock Pure discrete spectrum dynamical system and periodic tiling
  associated with a substitution.
\newblock {\em Ann. Inst. Fourier}, 54(2):288--299, 2004.

\bibitem[Sin06]{Sing:06}
B.~Sing.
\newblock Iterated function systems in mixed {E}uclidean and {$ p$}-adic
  spaces.
\newblock In {\em Complexus mundi}, pages 267--276. World Sci. Publ.,
  Hackensack, NJ, 2006.

\bibitem[Soc90]{Socolar}
J.~E.~S. Socolar.
\newblock Weak matching rules for quasicrystals.
\newblock {\em Comm. Math. Phys.}, 129(3):599--619, 1990.

\bibitem[Sol97]{Sol97}
B.~Solomyak.
\newblock Dynamics of self-similar tilings.
\newblock {\em Ergod. Th. Dyn. Sys.}, 17:695--738, 1997.

\bibitem[Sol07]{Sol:07}
B.~Solomyak.
\newblock Eigenfunctions for substitution tiling systems.
\newblock In {\em Probability and number theory---{K}anazawa 2005}, volume~49
  of {\em Adv. Stud. Pure Math.}, pages 433--454. Math. Soc. Japan, Tokyo,
  2007.

\bibitem[ST09]{SieThus:09}
A.~Siegel and J.~M. Thuswaldner.
\newblock Topological properties of {R}auzy fractals.
\newblock {\em M\'{e}m. Soc. Math. Fr. (N.S.)}, (118):140, 2009.

\bibitem[Thu89]{Thurston}
W.~P. Thurston.
\newblock Groups, tilings and finite state automata.
\newblock Lectures notes distributed in conjunction with the {C}olloquium
  {S}eries, in {\it {AMS} {C}olloquium Lectures}, 1989.

\bibitem[Thu20]{thuswaldner2019boldsymbolsadic}
J.~M. Thuswaldner.
\newblock {$S$}-adic sequences: a bridge between dynamics, arithmetic, and
  geometry.
\newblock In {\em Substitution and tiling dynamics: introduction to
  self-inducing structures}, volume 2273 of {\em Lecture Notes in Math.}, pages
  97--191. Springer, Cham, 2020.

\bibitem[Vee65]{Veech-1965}
W.~A. Veech.
\newblock Almost automorphic functions on groups.
\newblock {\em Amer. J. Math.}, 87:719--751, 1965.

\bibitem[Vee69]{Veech-1969}
W.~A. Veech.
\newblock Properties of minimal functions on abelian groups.
\newblock {\em Amer. J. Math.}, 91:415--440, 1969.

\bibitem[VGG04]{verger}
J.-L. Verger~Gaugry and J.-P. Gazeau.
\newblock Geometric study of the beta-integers for a {P}erron number and
  mathematical quasicrystals.
\newblock {\em J. Th\'eor. Nombres Bordeaux}, 16:125--149, 2004.

\bibitem[VM01]{VM:01}
J.~Vidal and R.~Mosseri.
\newblock Generalized quasiperiodic {R}auzy tilings.
\newblock {\em J. Phys. A}, 34(18):3927--3938, 2001.

\bibitem[Yoc06]{yoccoz}
J.-C. Yoccoz.
\newblock Continued fraction algorithms for interval exchange maps: an
  introduction.
\newblock In {\em Frontiers in number theory, physics, and geometry. {I}},
  pages 401--435. Springer, Berlin, 2006.

\end{thebibliography}

\end{document}